\documentclass[12pt,a4paper,oneside,reqno]{article}

\usepackage[utf8]{inputenc}
\usepackage[english]{babel}
\usepackage[left=1in,right=1in,top=1in,bottom=1in]{geometry}
\usepackage{amsfonts,amssymb,amsmath,amsthm}
\usepackage[symbol]{footmisc}
\usepackage{longtable,array}
\newcolumntype{L}[1]{>{\raggedright\arraybackslash$}p{#1}<{$}}
\usepackage{enumerate}
\usepackage{cite}
\usepackage{xcolor}
\usepackage[normalem]{ulem}
\usepackage{cancel}
\usepackage{tikz}
\usepackage{adjustbox}
\usetikzlibrary{arrows}
\usepackage{hyperref}
\hypersetup{hidelinks}
\usepackage{fancyhdr}
\theoremstyle{plain}
\newtheorem{theorem}{Theorem}
\newtheorem{lemma}[theorem]{Lemma}

\newtheorem{corollary}[theorem]{Corollary}

\numberwithin{theorem}{section}
\numberwithin{equation}{section}

\IfFileExists{fouriernc.sty}{\usepackage{fouriernc}}{}

\newcommand{\doi}[1]{\href{https://doi.org/#1}{doi:#1}}
\newcommand{\arxiv}[1]{\href{https://arxiv.org/abs/#1}{arXiv:#1}}
\newcommand{\arxivsearch}[1]{\href{https://arxiv.org/search/?query=#1}{arXiv}}
\tikzset{newedge/.style={}}

\begin{document}

\noindent\fbox{\parbox{\dimexpr\textwidth-2\fboxsep-2\fboxrule\relax}{\small
\textbf{Corrected version.} This is a corrected version of the paper
I.~Kaygorodov, Yu.~Volkov, \emph{The variety of $2$-dimensional algebras over an algebraically
closed field}, Canad. J. Math. 71 (2019), no.~4, 819--842, \doi{10.4153/S0008414X18000056}.}}

\bigskip

\noindent{\Large
The variety of $2$-dimensional algebras over an algebraically closed field}

\bigskip

\begin{center}
{\bf
Ivan Kaygorodov \& Yury Volkov}
\end{center}

\noindent {\bf Abstract.}
{\it The work is devoted to the variety of $2$-dimensional algebras over an algebraically closed field.
Firstly, we classify such algebras modulo isomorphism.
Then we describe the degenerations and the closures of certain algebra series in the variety of
$2$-dimensional algebras. Finally, we apply our results to obtain analogous descriptions for
the subvarieties of flexible and bicommutative algebras. In particular, we describe rigid algebras
and irreducible components for these subvarieties.}

\bigskip


\bigskip


\bigskip

\tableofcontents

\section*{Introduction}
\addcontentsline{toc}{section}{Introduction}

In this paper, an algebra is simply a vector space over a field with a bilinear binary operation that doesn't have to be associative.
Algebras of a fixed dimension form a variety with a natural action of a general linear group. Orbits under this action correspond to isomorphism classes of algebras. 
There are many classifications up to isomorphism for varieties of algebras of some fixed dimension satisfying some polynomial identities.
For example, there exist such classifications of
$3$-dimensional Novikov algebras \cite{bai},
$4$-dimensional Leibniz algebras \cite{Abror13},
$6$-dimensional Lie algebras \cite{SW} and many others.

In this paper we classify all $2$-dimensional algebras over an algebraically closed field up to isomorphism. It is not the first work devoted to this problem, classifications of different types were made in \cite{goze,GR11,mir00}, but 
all of them are not convenient for our main goal, the geometric description of the algebraic variety of $2$-dimensional algebras.
One of the advantages of our paper  is that our approach deals uniformly with all possible characteristics while the authors of \cite{mir00} don't consider the characteristics $2$ and $3$ and the authors of \cite{goze} consider only the two elements field in the characteristic $2$.
The authors of \cite{mir00} in fact don't give an explicit classification of $2$-dimensional algebras up to isomorphism because they have other purposes. They describe the moduli space by proving that $2$-dimensional algebras can be divided into parts that can be naturally included into projective spaces of different dimensions. The authors pretend that the classification up to isomorphism is easy and could be extracted from their proofs. The classification is really not very difficult and we believe that one can extract it from \cite{mir00} after reading the paper, taking parts of the classification from different places and taking in account carefully all the details while for us it was easier to produce this classification from scratch. The authors of \cite{goze} have produced a full classification. One of the problems is that this classification is outstretched through the whole paper and is mixed with other formulas. To collect all the parts of the classification from  \cite{goze} in one place and find all the additional conditions for these parts one has to fulfill a tedious work. At the same time, \cite{goze} contains some inaccuracies. For example, the series $\mu_{10}$ parametrized by two scalars has to be divided into two series parametrized by one scalar, the series $\mu_{11}$ admits nontrivial isomorphisms, and in the case of a commutative $2$-dimensional algebra with one idempotent $e$ it may be impossible to find $f$ linear independent with $e$ such that $f^2$ and $e$ are linearly dependent. The paper \cite{GR11} is very nice and gives the full classification of $2$-dimensional algebras over any field. Unfortunately, the answer is not given in terms of multiplication tables. The translation of this answer to the language of multiplication tables as well as its direct usage for the description of orbit closures is very difficult and it seems to be easier to produce a new appropriate classification. Also the consideration of arbitrary fields complicates the result and the extraction of the answer for an algebraically closed field becomes tedious.
For these reasons, we give a classification that is valid over an algebraically closed field of arbitrary characteristic.
In the same part of the paper, we also describe the automorphism groups for all algebras under consideration.

In the main part of our paper we develop the geometry of the variety of $2$-dimensional algebras. Namely, we describe the closures of orbits of some sets with respect to the Zariski topology.
Firstly, we describe all possible degenerations, i.e. closures of orbits of one point sets.
Degenerations are an interesting subject, which was studied in various papers (see, for example, \cite{B99,B05,M79,M80,BC99,S90,GRH,GRH2,GRH3,BB09,ikv17,BB14,laur03,kppv, kpv,gorb91,gorb93,gorb98}).
One of the problems in this direction is to describe all degenerations in a variety of algebras of some fixed dimension satisfying some set of identities.
For example, this problem was solved
for $2$-dimensional pre-Lie algebras in \cite{BB09},  
for $3$-dimensional Novikov algebras in \cite{BB14}, 
for $4$-dimensional Lie algebras in \cite{BC99}, 
for $4$-dimensional Zinbiel and nilpotent Leibniz algebras in \cite{kppv},
for nilpotent $5$- and $6$-dimensional Lie algebras in \cite{S90,GRH}, 
and for nilpotent $5$- and $6$-dimensional Malcev algebras in \cite{kpv}. As an application of our results, one can easily recover the results of \cite{BB09}.

Another interesting notion concerning degenerations is the so-called level of an algebra defined in the end of Section 5. The algebras of the first level and the associative, Lie and Jordan algebras of the second level are classified in \cite{khud13,khud15}.
In the papers \cite{gorb91,gorb93,gorb98}, the author defined the notion of an infinite level and described all anticommutative algebras that have an infinite level not greater than $3$. This notion is much easier in the sense that the infinite level of an algebra can be easily expressed in terms of the usual level. Algebras of low dimension play a special role in problems of such type, because they have small levels. The complete description of degenerations obtained in this work allows to compute the level for all $2$-dimensional algebras.

The next result of this paper is the description of orbit closures of certain series that appear in the classification up to isomorphism. Let $T$ be some subvariety of the variety of $n$-dimensional algebras closed under the action of the general linear group. An $n$-dimensional algebra from $T$ is called rigid if its orbit is an open subset of $T$.
Another important characteristic of a variety is its partition into irreducible components. The notion of a rigid algebra is closely related to this characteristic, because orbit closures of such algebras form irreducible components.
For example, irreducible components and rigid algebras were classified for low dimensional associative (see \cite{M79,M80}) and Jordan (see \cite{KE14}) algebras.
Since the variety of $2$-dimensional algebras is simply ${\bf k}^8$, it is clear that there is only one irreducible component and there are no rigid algebras in it.
Thus, this problem is not relevant for the variety of all $2$-dimensional algebras itself. Nevertheless, it is relevant for subvarieties.
In the last part we apply our results about the variety of all $2$-dimensional algebras to its subvarieties consisting of flexible and bicommutative algebras. We describe all degenerations and closures of orbits in these varieties.
In particular, we classify the irreducible components and rigid algebras. Our results allow to get such descriptions and classifications for varieties of $2$-dimensional algebras defined by any identities without any problems.

Let us give a resume of our motivations.
The problems considered in this paper are classical and their solution is interesting itself. In general the classification of all $n$-dimensional algebras is a wild problem and it is interesting to get the solution in particular cases where it is still possible. Our main motivation was the classification of all the algebras of the second level that we produce in \cite{kayvo} using the results of this paper. In fact, there are reasons to guess that our classification will allow to classify also algebras of the third and the fourth levels. Thus, our results are important for the classification of algebras of small levels and constitute a necessary part of it. Another application that we have in mind is the geometric description of subvarieties of the variety of $2$-dimensional algebras. There are some works (for example, \cite{BB09}) devoted to this problem and our work gives a powerful tool to solve it in all certain cases.
Our results can be applied whenever the natural action of $GL({\bf k}^2)$ on $({\bf k}^2)^* \otimes ({\bf k}^2)^* \otimes {\bf k}^2$ appears and we expect that they will have other applications, for example, in the theory of algebras with polynomial identities or in the geometric representation theory. Even when one deals with $n$-dimensional algebras for $n>2$ it may be useful to consider $2$-dimensional subalgebras, and our results could be applied in this case. For example, the classification of $n$-dimensional algebras with an $(n-2)$-dimensional annihilator is fulfilled using our classification in  \cite{bi2}.


\section{Definitions and notation}

Throughout  the paper we fix an algebraically closed field ${\bf k}$, a $2$-dimensional ${\bf k}$-linear vector space $V$ and a basis $\{e_1,e_2\}$ of $V$.
All spaces in this paper are considered over ${\bf k}$, and we write simply $dim$, $Hom$ and $\otimes$ instead of $dim_{{\bf k}}$, $Hom_{{\bf k}}$ and $\otimes_{{\bf k}}$. An algebra $A$ is a set with a structure of a vector space and a binary operation that induces a bilinear map from $A\times A$ to $A$.

Since this paper is devoted to $2$-dimensional algebras, we give all definitions and notation only for this case, though everything in this section can be rewritten for any dimension.

The set $\mathcal{A}_2:=Hom(V \otimes V,V) \cong V^* \otimes V^* \otimes V$ is a vector space of dimension $8$. This space has a structure of the affine variety ${\bf k}^8.$ Indeed, any $\mu\in \mathcal{A}_2$ is determined by $8$ structure constants $c_{ij}^k\in{\bf k}$ ($i,j,k=1,2$) such that
$\mu(e_i\otimes e_j)=c_{ij}^1e_1+c_{ij}^2e_2$. A subset of $\mathcal{A}_2$ is {\it Zariski-closed} if it can be defined by a set of polynomial equations in the variables $c_{ij}^k$.

The general linear group $GL(V)$ acts on $\mathcal{A}_2$ by conjugations:
$$ (g * \mu )(x\otimes y) = g\mu(g^{-1}x\otimes g^{-1}y)$$ 
for $x,y\in V$, $\mu\in \mathcal{A}_2$ and $g\in GL(V)$.
Thus, $\mathcal{A}_2$ is decomposed into $GL(V)$-orbits that correspond to the isomorphism classes of $2$-dimensional algebras. 
The classification of $2$-dimensional algebras up to isomorphism is equivalent to the classification of $GL(V)$-orbits.

Let $O(\mu)$ denote the orbit of $\mu\in\mathcal{A}_2$ under the action of $GL(V)$ and $\overline{O(\mu)}$ denote the Zariski closure of $O(\mu)$.
Let $A$ and $B$ be two $2$-dimensional algebras and $\mu,\lambda \in \mathcal{A}_2$ represent $A$ and $B$ respectively.
We say that $A$ degenerates to $B$ and write $A\to B$ if $\lambda\in\overline{O(\mu)}$.
Note that in this case we have $\overline{O(\lambda)}\subset\overline{O(\mu)}$. Hence, the definition of a degeneration doesn't depend on the choice of $\mu$ and $\lambda$. If $A\not\cong B$, then the assertion $A\to B$ is called a {\it proper degeneration}. We write $A\not\to B$ if $\lambda\not\in\overline{O(\mu)}$.
Let now $A(*):=\{A(\alpha)\}_{\alpha\in I}$ be a set of $2$-dimensional algebras and $\mu_\alpha\in\mathcal{A}_2$ represent $A(\alpha)$ for $\alpha\in I$.
If $\lambda\in\overline{\{O(\mu_\alpha)\}_{\alpha\in I}}$, then we write $A(*)\to B$ and say that $A(*)$ degenerates to $B$. In the opposite case we write $A(*)\not\to B$.

Let $A(*)$, $B$, $\mu_\alpha$ ($\alpha\in I$) and $\lambda$ be as above. Let $c_{ij}^k$ ($i,j,k=1,2$) be the structure constants of $\lambda$ in the basis $e_1,e_2$. If we construct maps $a_i^j:{\bf k}^*\to {\bf k}$ ($i,j=1,2$) and $f: {\bf k}^* \to I$ such that $a_1^1(t)e_1+a_1^2(t)e_2$ and $a_2^1(t)e_1+a_2^2(t)e_2$ form a basis of $V$ for any  $t\in{\bf k}^*$, and the structure constants of $\mu_{f(t)}$ in this basis are  polynomials $c_{ij}^k(t)\in{\bf k}[t]$ such  that $c_{ij}^k(0)=c_{ij}^k$, then $A(*)\to B$.
Indeed, if there is some closed subset $\mathcal{R}$ containing $O(\mu_{\alpha})$ for all $\alpha\in I$, then it contains, in particular, $O(\mu_{f(t)})$ for all $t\in{\bf k}^*$, and hence the element $\lambda_t$ of $\mathcal{A}_2$ with structure constants $c_{ij}^k(t)$ belongs to $\mathcal{R}$ for any $t\in{\bf k}^*$.
Note that the assertion $\lambda_t\in \mathcal{R}$ is equivalent to the annihilation of some set of polynomials in one variable in the point $t$. But if this set of polynomials vanishes for all $t\in{\bf k}^*$, then each of these polynomials has infinitely many roots, and hence it equals zero.
Thus, $t=0$ annihilates all the required polynomials too, i.e. $\lambda=\lambda_0\in \mathcal{R}$.
We will call  $(a_1^1(t)e_1+a_1^2(t)e_2,a_2^1(t)e_1+a_2^2(t)e_2)$ and $f(t)$ a {\it parametrized basis} and a {\it parametrized index} for $A(*)\to B$ respectively. The case of degeneration between two algebras corresponds to the case $|I|=1$. In this case we need only a parametrized basis, because $f(t)$ is the unique element of $I$ for any $t\in {\bf k}^*$.

We take the ideas for proving non-degenerations from \cite{S90}. Let $Q$ be a set of polynomial  equations in the variables $x_{i,j}^k$ ($i,j,k=1,2$). Suppose that $Q$ satisfies the following property: 
if $x_{i,j}^k=c_{ij}^k$ is a solution to all equations in $Q$, then  also $x_{i,j}^k=\tilde c_{ij}^k$ is a solution to all equations in $Q$ too in the following cases:
\begin{enumerate}
    \item there are $\alpha_1,\alpha_2\in{\bf k}^*$ such that $\tilde c_{ij}^k=\frac{\alpha_i\alpha_j}{\alpha_k}c_{ij}^k$;
    \item there is $\alpha\in{\bf k}$ such that
\begin{multline*}
\tilde c_{11}^1=c_{11}^1+\alpha(c_{12}^1+c_{21}^1)+\alpha^2c_{22}^1,\,\tilde c_{21}^1=c_{21}^1+\alpha c_{22}^1,\,\tilde c_{12}^1=c_{12}^1+\alpha c_{22}^1,\,\tilde c_{22}^1=c_{22}^1,\,\\
\tilde c_{11}^2=c_{11}^2+\alpha(c_{12}^2+c_{21}^2-c_{11}^1)+\alpha^2(c_{22}^2-c_{12}^1-c_{21}^1)-\alpha^3c_{22}^1,\\
\tilde c_{21}^2=c_{21}^2+\alpha(c_{22}^2-c_{21}^1)-\alpha^2c_{22}^1,\,
\tilde c_{12}^2=c_{12}^2+\alpha(c_{22}^2-c_{12}^1)-\alpha^2c_{22}^1,\,\tilde c_{22}^2=c_{22}^2-\alpha c_{22}^1.
\end{multline*}
\end{enumerate}
Let $\mathcal{R}\subset\mathcal{A}_2$ be a set of all algebra structures whose structure constants satisfy all equations in $Q$. We will call such a set $\mathcal{R}$ a {\it closed upper invariant set}. Let $\{A(\alpha)\}_{\alpha\in I}$ be a set of $2$-dimensional algebras such that $A(\alpha)$ can be represented by a structure from $\mathcal{R}$ for any $\alpha\in I$. Let $B$ be a $2$-dimensional algebra represented by the structure $\lambda\in\mathcal{A}_2$. If $O(\lambda)\cap\mathcal{R}=\varnothing$, then $A(*)\not\to B$. In this case we call $\mathcal{R}$ a {\it separating set} for $A(*)\not\to B$.

Let us recall two more tools for proving degenerations and non-degenerations. Firstly, if $A\to B$ is a proper degeneration, then ${\rm dim}\,{\rm Aut}(A)<{\rm dim}\,{\rm Aut}(B)$. Note that if $A(*)\to B$, then either
${\rm dim}\,{\rm Aut}(A(\alpha))={\rm dim}\,{\rm Aut}(B)$ for infinitely many $\alpha\in I$ or ${\rm dim}\,{\rm Aut}(A(\alpha))<{\rm dim}\,{\rm Aut}(B)$ for some $\alpha\in I$, but
it is possible that ${\rm dim}\,{\rm Aut}(A(\alpha))\ge {\rm dim}\,{\rm Aut}(B)$ for all $\alpha\in I$. Note also that ${\rm dim}\,{\rm Aut}(A)\mathrel{\le}{\rm dim}\,{\rm Der}(A)$, and the equality holds if $char{\bf k}=0$.
Secondly, if $A\to C$ and $C\to B$ then $A\to B$.
If there is no $C$ such that $A\to C$ and $C\to B$ are proper degenerations, then the assertion $A\to B$ is called a {\it primary degeneration}. If there are no $C$ and $D$ such that $C\to A$, $B\to D$, $C\not\to D$ and one of the assertions $C\to A$ and $B\to D$ is a proper degeneration,  then the assertion $A \not\to B$ is called a {\it primary non-degeneration}.
It suffices to prove only primary degenerations and non-degenerations to describe degenerations in the variety under consideration. Note also that any algebra degenerates to the algebra with zero multiplication.	

\section{Algebraic classification}

 The first of our aims is to classify all $2$-dimensional algebras over ${\bf k}$ modulo isomorphism. Our classification is based on the following lemma.
 
 \begin{lemma}\label{idem}
Let $A$ be a $2$-dimensional algebra. Then there exists a non-zero element $x\in A$ such that $x$ and $x^2$ are linearly dependent.
\end{lemma}
\begin{proof} The required assertion is equivalent to the existence of a $1$-dimensional subalgebra in $A$. Then the lemma follows from the discussion right after \cite[Proposition 1]{mir00}.

\end{proof}

Note that if $x\in A$ and $x^2$ are linearly dependent, then either $x^2=0$ or $x=\alpha e$ for some $\alpha\in{\bf k}^*$ and some $e\in A$ such that $e^2=e$. If $x^2=0$, then $x$ is called a {\it $2$-nil element}. An element $e$ such that $e^2=e$ is called an {\it idempotent}.

\begin{corollary}\label{classes}
Any $2$-dimensional ${\bf k}$-algebra belongs to one of the following disjoint classes:
\begin{enumerate}
\item[${\bf A}.$] algebras that don't have nonzero idempotents and have a unique $1$-dimensional subspace of $2$-nil elements;
\item[${\bf B}.$] algebras that don't have nonzero idempotents and have two linearly independent $2$-nil elements;
\item[${\bf C}.$] algebras that have a unique nonzero idempotent and don't have nonzero $2$-nil elements;
\item[${\bf D}.$] algebras that have a unique nonzero idempotent and a nonzero $2$-nil element;
\item[${\bf E}.$] algebras that have two different nonzero idempotents.
\end{enumerate}
\end{corollary}
\begin{proof} The fact that the classes are disjoint is obvious. The fact that any $2$-dimensional algebra belongs to one of the classes follows easily from Lemma \ref{idem} and the remark after it.
\end{proof}

To give the classification of $2$-dimensional algebras we have to introduce some notation.
Let us consider the action of the cyclic group $C_2=\langle \rho\mid \rho^2\rangle$ on ${\bf k}$ defined by the equality ${}^{\rho}\alpha=-\alpha$ for $\alpha\in{\bf k}$.
Let us fix some set of representatives of orbits under this action and denote it by ${\bf k_{\ge 0}}$. For example, if ${\bf k}=\mathbb{C}$, then one can take $\mathbb{C}_{\ge 0}=\{\alpha\in\mathbb{C}\mid Re(\alpha)>0\}\cup\{\alpha\in\mathbb{C}\mid Re(\alpha)=0,Im(\alpha)\ge 0\}$.

Let us also consider the action of $C_2$ on ${\bf k}^2$ defined by the equality ${}^{\rho}(\alpha,\beta)=(1-\alpha+\beta,\beta)$ for $(\alpha,\beta)\in{\bf k}^2$.
Let us fix some set of representatives of orbits under this action and denote it by $\mathcal{U}$. Let us also define $\mathcal{T}=\{(\alpha,\beta)\in{\bf k}^2\mid \alpha+\beta=1\}$.

Given $(\alpha,\beta,\gamma,\delta)\in{\bf k}^4$, we define $\mathcal{D}(\alpha,\beta,\gamma,\delta)=(\alpha+\gamma)(\beta+\delta)-1$.
We define $\mathcal{C}_1(\alpha,\beta,\gamma,\delta)=(\beta,\delta)$, $\mathcal{C}_2(\alpha,\beta,\gamma,\delta)=(\gamma,\alpha)$, and  $\mathcal{C}_3(\alpha,\beta,\gamma,\delta)=\left(\frac{\beta\gamma-(\alpha-1)(\delta-1)}{\mathcal{D}(\alpha,\beta,\gamma,\delta)},\frac{\alpha\delta-(\beta-1)(\gamma-1)}{\mathcal{D}(\alpha,\beta,\gamma,\delta)}\right)$ for $(\alpha,\beta,\gamma,\delta)$ such that $\mathcal{D}(\alpha,\beta,\gamma,\delta)\not=0$. Let us consider the set \begin{center}
$X=\left\{\big(\mathcal{C}_1(\Gamma),\mathcal{C}_2(\Gamma),\mathcal{C}_3(\Gamma)\big)\mid \Gamma\in{\bf k}^4, \mathcal{D}(\Gamma)\not=0,\mathcal{C}_1(\Gamma),\mathcal{C}_2(\Gamma)\not\in \mathcal{T}\right\}\subset ({\bf k}^2)^3.$
\end{center}
One can show that the symmetric group $S_3$ acts on $X$ by the equality  $${}^{\sigma}\big(\mathcal{C}_1(\Gamma),\mathcal{C}_2(\Gamma),\mathcal{C}_3(\Gamma)\big)=\big(\mathcal{C}_{\sigma^{-1}(1)}(\Gamma),\mathcal{C}_{\sigma^{-1}(2)}(\Gamma),\mathcal{C}_{\sigma^{-1}(3)}(\Gamma)\big) \mbox{ for }\sigma\in S_3.$$
Indeed, suppose that $\big(\mathcal{C}_1(\Gamma),\mathcal{C}_2(\Gamma),\mathcal{C}_3(\Gamma)\big)\in X$ for some $\Gamma=(\alpha,\beta,\gamma,\delta)$. We need to show that ${}^{\sigma}\big(\mathcal{C}_1(\Gamma),\mathcal{C}_2(\Gamma),\mathcal{C}_3(\Gamma)\big)\in X$ for any $\sigma\in S_3$. We will check this for $\sigma$ interchanging $1$ and $3$, the other verifications are analogous. First, $\mathcal{C}_3(\Gamma)\in\mathcal{T}$ is equivalent to the equality $\alpha+\beta+\delta+\gamma-2=\mathcal{D}(\Gamma)$ that can be rewritten in the form $(\alpha+\gamma-1)(\beta+\delta-1)=0$. It is clear that the last equality is not valid. Let us introduce
$$\Gamma'=\left(\alpha,\frac{\beta\gamma-(\alpha-1)(\delta-1)}{\mathcal{D}(\Gamma)},\gamma,\frac{\alpha\delta-(\beta-1)(\gamma-1)}{\mathcal{D}(\Gamma)}\right).$$
It remains to check that $\mathcal{D}(\Gamma')\not=0$ and $\mathcal{C}_3(\Gamma')=(\beta,\delta)$. The equality $\mathcal{D}(\Gamma')=0$ is equivalent to the equality $(\alpha+\beta+\delta+\gamma-2)(\alpha+\gamma)=\mathcal{D}(\Gamma)$ that can be rewritten in the form $(\alpha+\gamma-1)^2=0$.
Hence, we get $\mathcal{D}(\Gamma')\not=0$. To prove that $\mathcal{C}_3(\Gamma')=(\beta,\delta)$ we have to verify two equalities. We will consider only the first equality, the second one is analogous. Thus, it remains to show that
\begin{multline*}
\frac{\beta\gamma-(\alpha-1)(\delta-1)}{\mathcal{D}(\Gamma)}\gamma-(\alpha-1)\left(\frac{\alpha\delta-(\beta-1)(\gamma-1)}{\mathcal{D}(\Gamma)}-1\right)=  \beta\left(\frac{(\alpha+\beta+\gamma+\delta-2)(\alpha+\gamma)}{\mathcal{D}(\Gamma)}-1\right).
\end{multline*}
Multiplying by $\mathcal{D}(\Gamma)$ and reducing all the equal terms, one sees that the last equality is valid.
Note that there exists a set of representatives of orbits $\mathcal{\tilde V}$ under the action of $S_3$ on $X$ such that if $(\mathcal{C}_1,\mathcal{C}_2,\mathcal{C}_3)\in \mathcal{\tilde V}$ and $\mathcal{C}_1\not=\mathcal{C}_2$, then $\mathcal{C}_3\not=\mathcal{C}_1,\mathcal{C}_2$. Let us fix such $\mathcal{\tilde V}$ and define
$$
\mathcal{V}=\{\Gamma\in{\bf k}^4\mid \mathcal{D}(\Gamma)\not=0; \mathcal{C}_1(\Gamma),\mathcal{C}_2(\Gamma)\not\in \mathcal{T}, \big(\mathcal{C}_1(\Gamma),\mathcal{C}_2(\Gamma),\mathcal{C}_3(\Gamma)\big)\in\mathcal{\tilde V}\}.
$$
For $\Gamma\in\mathcal{V}$, we also define $\mathcal{C}(\Gamma)=\{\mathcal{C}_1(\Gamma),\mathcal{C}_2(\Gamma),\mathcal{C}_3(\Gamma)\}\subset{\bf k}^2$.

Let us consider the action of the cyclic group $C_2$ on ${\bf k}^*\setminus \{1\}$ defined by the equality ${}^{\rho}\alpha=\alpha^{-1}$ for $\alpha\in{\bf k}^*\setminus \{1\}$.
Let us fix some set of representatives of orbits under this action and denote it by ${\bf k_{>1}^*}$. For example, if ${\bf k}=\mathbb{C}$, then one can take 
\begin{center}$\mathbb{C}_{>1}^*=\{\alpha\in\mathbb{C}^*\mid |\alpha|>1\}\cup\{\alpha\in\mathbb{C}^*\mid |\alpha|=1,0<arg(\alpha)\le \pi\}$.\end{center} For $(\alpha,\beta,\gamma)\in {\bf k}^2\times{\bf k}^*_{>1}$ we define $$\mathcal{C}(\alpha,\beta,\gamma)=\left\{\big(\alpha\gamma,(1-\alpha)\gamma\big),\left(\frac{\beta}{\gamma},\frac{1-\beta}{\gamma}\right)\right\}\subset{\bf k}^2.$$
Finally, let us consider the action of $C_2$ on ${\bf k}^2$ defined by the equality ${}^{\rho}(\alpha,\beta)=(\beta,\alpha)$, fix some set of representatives of orbits under this action, denote it by $\mathcal{S}$, and define
$$\mathcal{W}=\{(\alpha,\beta,\gamma)\in{\bf k}^2\times{\bf k}^*_{>1}\mid \gamma\not=-1\mbox{ or }(\alpha,\beta)\in\mathcal{S}\}.$$
Note that $-1\in{\bf k}^*_{>1}$ if $char{\bf k}\not=2$, because $-1$ is a fixed point of the action ${}^{\rho}\alpha=\alpha^{-1}$, and that $\mathcal{W}={\bf k}^2\times{\bf k}^*_{>1}$ if $char{\bf k}=2$.

Let $\mathcal{F}\subset \mathcal{A}_2$ be the set formed by the algebra structures on the vector space $V$ listed in Table 1.
This section is devoted to the proof of the following theorem that gives a classification of $2$-dimensional algebras over ${\bf k}$ up to isomorphism.

\begin{theorem}\label{alg}
Any non-trivial $2$-dimensional ${\bf k}$-algebra can be represented by a unique structure from $\mathcal{F}$.
\end{theorem}

In other words, Theorem \ref{alg} states that $\mathcal{A}_2=\bigcup\limits_{\mu\in\mathcal{F}}O(\mu)\cup\{{\bf k^2}\}$ and that, if $\mu,\lambda\in\mathcal{F}$ are different structures, then $O(\mu)\cap O(\lambda)=\varnothing$.
Whenever an algebra named $A$ appears in this section, we suppose that it is represented by some structure from $\mathcal{A}_2$ with structure constants $c_{ij}^k$ ($i,j,k=1,2$).
 According to Corollary \ref{classes}, it suffices to consider each of the classes ${\bf A}$--${\bf E}$ separately. It is not difficult to show that the letter in the name of an algebra from $\mathcal{F}$ corresponds to its class in each case. This will follow also from our proofs.

\begin{lemma}\label{Aclass}
If $A$ belongs to the class ${\bf A}$, then it can be represented by a unique structure from the set
\begin{equation}\label{Aset}
\{{\bf A}_1(\alpha)\}_{\alpha\in{\bf k}}\cup \{{\bf A}_2\}\cup\{{\bf A}_3\}\cup\{{\bf A}_4(\alpha)\}_{\alpha\in{\bf k_{\ge 0}}}.
\end{equation}
\end{lemma}
\begin{proof}
Let us represent the algebra $A$ by a structure such that $e_2e_2=0$. It is easy to see that $A$ belongs to the class ${\bf A}$ iff $x_t=e_1+te_2$ and $x_t^2$ are linearly independent for any $t\in{\bf k}$. Since 
$$x_t^2=(c_{11}^1+(c_{12}^1+c_{21}^1)t)e_1+(c_{11}^2+(c_{12}^2+c_{21}^2)t)e_2,$$
$x_t$ and $x_t^2$ are linearly independent iff
$$
0\not=\left|\begin{array}{cc}
c_{11}^1+(c_{12}^1+c_{21}^1)t&c_{11}^2+(c_{12}^2+c_{21}^2)t\\
1&t
\end{array}\right|=(c_{12}^1+c_{21}^1)t^2+(c_{11}^1-c_{12}^2-c_{21}^2)t-c_{11}^2.
$$
Since by our assumption $x_t$ and $x_t^2$ are linearly independent for any $t\in{\bf k}$, we have $c_{12}^1+c_{21}^1 =0$,
$c_{11}^1=c_{12}^2+c_{21}^2$, and $c_{11}^2\not=0$.

Now we have four cases:
\begin{itemize}
\item $c_{12}^1=0$, $c_{11}^1\neq 0$. Considering the basis $\frac{e_1}{c_{11}^1}$, $\frac{c_{11}^2e_2}{( c_{11}^1)^2} $ of $V$, one can check that $A$ can be represented by ${\bf A}_{1}\left(\frac{c_{12}^2}{c_{11}^1}\right)$.

\item $c_{12}^1=0$, $c_{12}^2=-c_{21}^2\not=0$. Considering the basis $\frac{e_1}{c_{12}^2}$, $\frac{c_{11}^2e_2}{( c_{12}^2)^2}$ of $V$, one can check that $A$ can be represented by ${\bf A}_{2}$.

\item $c_{12}^1=c_{12}^2=c_{21}^2=0$. Considering the basis $e_1$, $c_{11}^2e_2$ of $V$, one can check that $A$ can be represented by ${\bf A}_{3}$.

\item $c_{12}^1\neq 0$. Let $a\in{\bf k}^*$ be such that $c_{11}^2c_{12}^1a^2=1$ and $c_{11}^1a\in{\bf k_{\ge 0}}$.
Considering the basis $a\left(e_1-\frac{c_{21}^2}{c_{12}^1}e_2\right)$, $\frac{e_2}{c_{12}^1}$ of $V$, one can check that $A$ can be represented by ${\bf A}_4(c_{11}^1a)$.
\end{itemize}

It remains to prove that any two different structures from the set \eqref{Aset} represent non-isomorphic algebras. Firstly, note that ${\rm dim} ({\bf A}_2)^2={\rm dim} ({\bf A}_3)^2=1$ while ${\rm dim} ({\bf A}_1(\alpha))^2={\rm dim} ({\bf A}_4(\alpha))^2=2$ for any $\alpha\in{\bf k}$. We have also ${\bf A}_2\not\cong {\bf A}_3$, because ${\bf A}_3$ has a nonzero annihilator. 

Suppose that $A$ is represented by the structure ${\bf A}_1(\alpha)$ for some $\alpha\in{\bf k}$. Then there exists $x\in A$ such that $x^2=0$, $xA+Ax\subset \langle x\rangle$, and $\alpha xy=(1-\alpha) yx$ for any $y\in A$. Such an element doesn't exist in ${\bf A}_4(\beta)$ for any $\beta\in{\bf k}$ and in ${\bf A}_1(\beta)$ for any $\beta\in{\bf k}\setminus\{\alpha\}$.

Suppose that $A$ is represented by the structure ${\bf A}_4(\alpha)$ for some $\alpha\in{\bf k_{\ge 0}}$.
Suppose that the structure constants of $A$ in the basis $E_1$, $E_2$ equal the structure constants of ${\bf A}_4(\beta)$ for some $\beta\in{\bf k_{\ge 0}}$. Since $E_2E_2=0$ and $E_2E_1=-E_1$, it is easy to see that $E_2=e_2$ and $E_1=ae_1$ for some $a\in{\bf k}^*$. Then we obtain from the equality $E_1E_1=\beta E_1+E_2$ that $a=\pm 1$ and $\beta=\pm \alpha$. Since $\alpha,\beta\in{\bf k_{\ge 0}}$, we have $\beta=\alpha$.
\end{proof}

\begin{lemma}\label{Bclass}
If $A$ belongs to the class ${\bf B}$, then it can be represented by a unique structure from the set
\begin{equation}\label{Bset}
\{{\bf B}_1(\alpha)\}_{\alpha\in{\bf k}}\cup\{{\bf B}_2(\alpha)\}_{\alpha\in{\bf k}}\cup \{{\bf B}_3\}.
\end{equation}
\end{lemma}
\begin{proof}
Let us represent the algebra $A$ by a structure such that $e_1e_1=e_2e_2=0$. For $s,t\in{\bf k}$, let us define $x_{s,t}=se_1+te_2$.
Suppose that there are $s,t\in{\bf k}^*$ such that $0=x_{s,t}^2=st(e_1e_2+e_2e_1)$. Then $e_1e_2+e_2e_1=0$ and $A$ is anticommutative. It is easy to see that any $2$-dimensional anticommutative algebra either has the trivial multiplication or can be represented by ${\bf B}_3$ (note that by our definition $A$ is anticommutative iff $x^2=0$ for any $x\in A$).

Suppose now that $x_{s,t}^2\not=0$ for any $s,t\in{\bf k}^*$. Since $A$ doesn't have idempotents, $x_{s,t}$ and $x_{s,t}^2$ are linearly independent for $s,t\in{\bf k}^*$.
It is easy to check that $x_{s,t}$ and $x_{s,t}^2$ are linearly dependent for $s=c_{12}^1+c_{21}^1$, $t=c_{12}^2+c_{21}^2$.
Hence, $c_{12}^1+c_{21}^1=0$ or $c_{12}^2+c_{21}^2=0$. Without loss of generality we may assume that $c_{12}^2+c_{21}^2=0$.
Since $A$ is not anticommutative, we have $c_{12}^1+c_{21}^1\not=0$ in this case.

If $c_{12}^2\not=0$, then, considering the basis $\frac{e_1}{c_{12}^2}$, $\frac{e_2}{c_{12}^1+c_{21}^1}$ of $V$, one can check that $A$ can be represented by
${\bf B}_{1}\left(\frac{c_{21}^1}{c_{12}^1+c_{21}^1}\right)$. If $c_{12}^2=0$, then, considering the basis $e_1$, $\frac{e_2}{c_{12}^1+c_{21}^1}$ of $V$, one can check that $A$ can be represented by ${\bf B}_{2}\left(\frac{c_{21}^1}{c_{12}^1+c_{21}^1}\right)$.

It remains to prove that any two different structures from the set \eqref{Bset} represent non-isomorphic algebras. Since ${\bf B}_3$ is anticommutative, it is not isomorphic to other algebras from \eqref{Bset}. Note also that ${\rm dim} ({\bf B}_1(\alpha))^2=2>1={\rm dim} ({\bf B}_2(\beta))^2$ for any $\alpha,\beta\in{\bf k}$.

Suppose that $A$ is represented by the structure ${\bf B}_i(\alpha)$ for some $\alpha\in{\bf k}$ and $i=1,2$.
Suppose that the structure constants of $A$ in the basis $E_1$, $E_2$ equal the structure constants of ${\bf B}_i(\beta)$ for some $\beta\in{\bf k}$.
Since $E_1E_1=E_2E_2=0$, we have either $E_1=ae_1$, $E_2=be_2$ or $E_1=ae_2$, $E_2=be_1$ for some $a,b\in{\bf k}^*$. Since $E_1E_2+E_2E_1=E_1$, we have
$E_1=ae_1$ and $E_2=e_2$. Then we get $\beta=\alpha$ from the equality $E_1E_2=(1-\beta)E_1+(2-i)E_2$.
\end{proof}

\begin{lemma}\label{Cclass}
If $A$ belongs to the class ${\bf C}$, then it can be represented by ${\bf C}(\alpha,\beta)$ for a unique pair $(\alpha,\beta)\in{\bf k}\times{\bf k_{\ge 0}}$.
\end{lemma}
\begin{proof}
Let us represent the algebra $A$ by a structure such that $e_2e_2=e_2$.  It is easy to see that $A$ belongs to the class ${\bf C}$ iff $x_t=e_1+te_2$ and $x_t^2$ are linearly independent for any $t\in{\bf k}$. Since 
$$x_t^2=(c_{11}^1+(c_{12}^1+c_{21}^1)t)e_1+(c_{11}^2+(c_{12}^2+c_{21}^2)t+t^2)e_2,$$
$x_t$ and $x_t^2$ are linearly independent iff
$$
0\not=\left|\begin{array}{cc}
c_{11}^1+(c_{12}^1+c_{21}^1)t&c_{11}^2+(c_{12}^2+c_{21}^2)t+t^2\\
1&t
\end{array}\right|=
(c_{12}^1+c_{21}^1-1)t^2+(c_{11}^1-c_{12}^2-c_{21}^2)t-c_{11}^2.
$$
Since by our assumption $x_t$ and $x_t^2$ are linearly independent for any $t\in{\bf k}$, we have 
\begin{center}$c_{12}^1+c_{21}^1 =1$,
$c_{11}^1=c_{12}^2+c_{21}^2$, and $c_{11}^2\not=0$.\end{center}

Let $a$ be such an element of ${\bf k}^*$ that $c_{11}^2a^2=1$ and $a(c_{12}^2-c_{11}^1c_{21}^1)\in{\bf k_{\ge 0}}$.
Considering the basis $a(e_1-c_{11}^1e_2)$, $e_2$ of $V$, one can check that $A$ can be represented by ${\bf C}\left(c_{21}^1,a(c_{12}^2-c_{11}^1c_{21}^1)\right)$.

Suppose that $A$ is represented by the structure ${\bf C}(\alpha,\beta)$ for some pair $(\alpha,\beta)\in{\bf k}\times{\bf k_{\ge 0}}$.
Suppose that the structure constants of $A$ in the basis $E_1$, $E_2$ equal the structure constants of ${\bf C}(\gamma,\delta)$ for some $(\gamma,\delta)\in{\bf k}\times{\bf k_{\ge 0}}$.
Since $E_2E_2=E_2$ and ${\bf C}(\alpha,\beta)$ has a unique idempotent, we have $E_2=e_2$. We get $E_1=\pm e_1$ from the equality $E_1E_1=E_2$.
Then $\gamma= \alpha$ and $\delta=\pm \beta$. Since $\beta,\delta\in{\bf k_{\ge 0}}$, we have $(\gamma,\delta)=(\alpha,\beta)$.
\end{proof}

\begin{lemma}\label{Dclass}
If $A$ belongs to the class ${\bf D}$, then it can be represented by a unique structure from the set
\begin{equation}\label{Dset}
\{{\bf D}_1(\alpha,\beta)\}_{(\alpha,\beta)\in \mathcal{U}}\cup\{{\bf D}_2(\alpha,\beta)\}_{(\alpha,\beta)\in{\bf k}^2\setminus \mathcal{T}}\cup\{{\bf D}_3(\alpha,\beta)\}_{(\alpha,\beta)\in{\bf k}^2\setminus \mathcal{T}}\cup\{{\bf D}_4(\alpha,\beta)\}_{(\alpha,\beta)\in{\bf k}^2}.
\end{equation}
\end{lemma}
\begin{proof}
Let us represent the algebra $A$ by a structure such that $e_1e_1=e_1$ and $e_2e_2=0$.

Let us consider the following cases:
\begin{itemize}
\item $c_{12}^1+c_{21}^1\not=0$. If $c_{12}^2+c_{21}^2\not\in\{0,1\}$, then one can check that $\frac{1}{c_{12}^2+c_{21}^2}\left(e_1+\frac{c_{12}^2+c_{21}^2-1}{c_{12}^1+c_{21}^1}e_2\right)$ is an idempotent that is not equal to $e_1$. Thus, $c_{12}^2+c_{21}^2\in\{0,1\}$. If $c_{12}^2+c_{21}^2=1$, then, considering the basis $e_1$, $\frac{e_2}{c_{12}^1+c_{21}^1}$ of $V$, one can check that $A$ can be represented by ${\bf D}_{4}\left(\frac{c_{21}^1}{c_{12}^1+c_{21}^1},c_{12}^2\right)$. Suppose now that $c_{12}^2+c_{21}^2=0$.
If $\left(\frac{c_{21}^1}{c_{12}^1+c_{21}^1},c_{12}^2\right)\in \mathcal{U}$, then, considering the basis $e_1$, $\frac{e_2}{c_{12}^1+c_{21}^1} $ of $V$, one can check that $A$ can be represented by ${\bf D}_{1}\left(\frac{c_{21}^1}{c_{12}^1+c_{21}^1},c_{12}^2\right)$.
If $\left(\frac{c_{21}^1}{c_{12}^1+c_{21}^1},c_{12}^2\right)\not\in \mathcal{U}$, then $\left(\frac{c_{12}^1}{c_{12}^1+c_{21}^1}+c_{12}^2,c_{12}^2\right)\in \mathcal{U}$ and, considering the basis $e_1$, $e_1-\frac{e_2}{c_{12}^1+c_{21}^1} $ of $V$, one can check that $A$ can be represented by ${\bf D}_{1}\left(\frac{c_{12}^1}{c_{12}^1+c_{21}^1}+c_{12}^2,c_{12}^2\right)$.

\item $c_{12}^1=-c_{21}^1\not=0$. Considering the basis $e_1$, $\frac{e_2}{c_{12}^1} $ of $V$, one can check that $A$ can be represented by ${\bf D}_{3}(c_{12}^2,c_{21}^2)$.
Since $e_1+e_2$ is not idempotent, $(c_{12}^2,c_{21}^2)\not\in \mathcal{T}$.

\item $c_{12}^1=c_{21}^1=0$. Then one can check that $A$ is represented by ${\bf D}_{2}(c_{12}^2,c_{21}^2)$.
Since $e_1+e_2$ is not idempotent, $(c_{12}^2,c_{21}^2)\not\in \mathcal{T}$.
\end{itemize}

It remains to prove that any two different structures from the set \eqref{Dset} represent non-isomorphic algebras.

Suppose that $A$ is represented by the structure ${\bf D}_1(\alpha,\beta)$ for some pair $(\alpha,\beta)\in \mathcal{U}$.
Note that $e_2e_2=(e_1-e_2)^2=0$ in ${\bf D}_1(\alpha,\beta)$ while the structures ${\bf D}_2(\gamma,\delta)$, ${\bf D}_3(\gamma,\delta)$ and ${\bf D}_4(\gamma,\delta)$ have a unique $1$-dimensional subspace of $2$-nil elements for any pair $(\gamma,\delta)\in{\bf k}^2$.
Suppose now that the structure constants of $A$ in the basis $E_1$, $E_2$ equal the structure constants of ${\bf D}_1(\gamma,\delta)$ for some pair $(\gamma,\delta)\in \mathcal{U}$. 
Since $E_1$ is an idempotent and ${\bf D}_1(\alpha,\beta)$ has a unique idempotent, we have $E_1=e_1$. Since $E_2E_2=0$, we have either $E_2=ae_2$ or $E_2=a(e_1-e_2)$ for some $a\in{\bf k}^*$. We obtain $a=1$ in both cases from the equality $E_1E_2+E_2E_1=E_1$. Then we have $\delta=\beta$ and either $\gamma=\alpha$ or $\gamma=1-\alpha+\beta$. Since
$(\alpha,\beta),(\gamma,\delta)\in \mathcal{U}$, we have $(\gamma,\delta)=(\alpha,\beta)$.

Suppose that $A$ is represented by the structure ${\bf D}_2(\alpha,\beta)$ for some pair $(\alpha,\beta)\in{\bf k}^2$. Note that $A$ has an element $x$ such that $x^2=0$ and $xA+Ax\subset\langle x\rangle$ while ${\bf D}_3(\gamma,\delta)$ and ${\bf D}_4(\gamma,\delta)$ don't have such an element for any pair $(\gamma,\delta)\in{\bf k}^2$ because any square zero element of ${\bf D}_3(\gamma,\delta)$ and of ${\bf D}_4(\gamma,\delta)$ is linearly dependent with $e_2$.
Suppose now that the structure constants of $A$ in the basis $E_1$, $E_2$ equal the structure constants of ${\bf D}_2(\gamma,\delta)$ for some pair $(\gamma,\delta)\in{\bf k}^2$. 
Since $E_1$ is an idempotent and ${\bf D}_2(\alpha,\beta)$ has a unique idempotent, we have $E_1=e_1$. Since $E_2E_2=0$, we have $E_2=ae_2$ for some $a\in{\bf k}^*$.
Then it is easy to see that $(\gamma,\delta)=(\alpha,\beta)$.

Finally, suppose that $A$ is represented by the structure ${\bf D}_3(\alpha,\beta)$ for some pair $(\alpha,\beta)\in{\bf k}^2$.
Suppose that the structure constants of $A$ in the basis $E_1$, $E_2$ equal the structure constants of ${\bf D}_3(\gamma,\delta)$ for some pair $(\gamma,\delta)\in{\bf k}^2$. 
Since $E_1$ is an idempotent and ${\bf D}_3(\alpha,\beta)$ has a unique idempotent, we have $E_1=e_1$. Since $E_2E_2=0$, we have $E_2=ae_2$ for some $a\in{\bf k}^*$.
Then it is easy to see that $a=1$ and $(\gamma,\delta)=(\alpha,\beta)$.

Note that ${\bf D}_3(\alpha,\beta)\not\cong{\bf D}_4(\gamma,\delta)$: for the unique nonzero idempotent $e$ and any nonzero square zero element $x$ we have $ex+xe\in\langle x\rangle$ in ${\bf D}_3(\alpha,\beta)$ (namely, $e_1e_2+e_2e_1=(\alpha+\beta)e_2$), while $e_1e_2+e_2e_1=e_1+e_2\not\in\langle e_2\rangle$ in ${\bf D}_4(\gamma,\delta)$.

Finally, suppose that $A$ is represented by the structure ${\bf D}_4(\alpha,\beta)$ and the structure constants of $A$ in the basis $E_1$, $E_2$ equal the structure constants of ${\bf D}_4(\gamma,\delta)$ for some pair $(\gamma,\delta)\in{\bf k}^2$. As above, $E_1=e_1$ and $E_2=ae_2$ for some $a\in{\bf k}^*$. The equality $E_1E_2+E_2E_1=E_1+E_2$ gives $a=1$, and hence $(\gamma,\delta)=(\alpha,\beta)$.
\end{proof}

As a consequence of the proofs of Lemmas \ref{Aclass}--\ref{Dclass} we can describe the automorphism groups of the algebras of the classes ${\bf A}$--${\bf D}$.

\begin{corollary}
$1.$ ${\rm Aut}({\bf A}_1(\alpha))\cong Aut({\bf A}_2)$ is isomorphic to the additive group of ${\bf k}$.

$2.$ ${\rm Aut}({\bf A}_3)$ is isomorphic to the subgroup of $GL_2({\bf k})$ formed by matrices of the form $\left(\begin{array}{cc}
a&0\\
b&a^2
\end{array}\right)$, where $a\in{\bf k}^*$ and $b\in{\bf k}$.

$3.$ ${\rm Aut}({\bf A}_4(\alpha))\cong C_2$ if $\alpha=0$ and $char{\bf k}\not=2$; ${\rm Aut}({\bf A}_4(\alpha))$ is trivial if either $\alpha\in{\bf k}^*$ or $\alpha=0$ and $char{\bf k}=2$.

$4.$ ${\rm Aut}({\bf B}_1(\alpha))$ is trivial; ${\rm Aut}({\bf B}_2(\alpha))\cong {\bf k}^*$.

$5.$ ${\rm Aut}({\bf B}_3)$ is isomorphic to the subgroup of $GL_2({\bf k})$ formed by matrices of the form $\left(\begin{array}{cc}
1&0\\
b&a
\end{array}\right)$, where $a\in{\bf k}^*$ and $b\in{\bf k}$.

$6.$ ${\rm Aut}({\bf C}(\alpha, \beta))\cong C_2$ if $\beta=0$ and $char{\bf k}\not=2$; ${\rm Aut}({\bf C}(\alpha,\beta))$ is trivial if either $\beta\in{\bf k}^*$ or $\beta=0$ and $char{\bf k}=2$.

$7.$ ${\rm Aut}({\bf D}_1(\alpha, \beta))\cong C_2$ if $\beta=2\alpha-1$ and ${\rm Aut}({\bf D}_1(\alpha,\beta))$ is trivial if $\beta\not=2\alpha-1$.

$8.$ ${\rm Aut}({\bf D}_2(\alpha, \beta))\cong {\bf k}^*$ and ${\rm Aut}({\bf D}_3(\alpha, \beta))$ is trivial if $\alpha+\beta\not=1$.

$9.$ ${\rm Aut}({\bf D}_4(\alpha, \beta))$ is trivial.

In particular, 
$${\rm dim}\,{\rm Aut}({\bf A}_4(\alpha))={\rm dim}\,{\rm Aut}({\bf B}_1(\alpha))={\rm dim}\,{\rm Aut}({\bf C}(\alpha,\beta))$$ $$={\rm dim}\,{\rm Aut}({\bf D}_1(\alpha,\beta))={\rm dim}\,{\rm Aut}({\bf D}_3(\alpha,\beta))={\rm dim}\,{\rm Aut}({\bf D}_4(\alpha,\beta))=0;$$
$${\rm dim}\,{\rm Aut}({\bf A}_1(\alpha))={\rm dim}\,{\rm Aut}({\bf A}_2)={\rm dim}\,{\rm Aut}({\bf B}_2(\alpha))={\rm dim}\,{\rm Aut}({\bf D}_2(\alpha, \beta))=1;$$ $${\rm dim}\,{\rm Aut}({\bf A}_3)={\rm dim}\,{\rm Aut}({\bf B}_3)=2.$$
\end{corollary}
\begin{proof}
1--3. Any structure of the class ${\bf A}$ has a unique subspace of $2$-nil elements generated by $e_2$. Thus, any automorphism of such an algebra sends $e_1$ and $e_2$ to $ae_1+be_2$ and $ce_2$ respectively, where $a,c\in{\bf k}^*$ and $b\in{\bf k}$. It is easy to check that $a=c=1$ for ${\bf A}_1(\alpha)$ and ${\bf A}_2$; $c=a^2$ for ${\bf A}_3$; $a=\pm 1$, $b=0$ and $c=1$ for ${\bf A}_4(0)$; and $a=c=1$, $b=0$ for ${\bf A}_4(\alpha)$ if $\alpha\not=0$.

4. It follows from the proof of Lemma \ref{Bclass} that any automorphism of the algebra ${\bf B}_i(\alpha)$, where $i\in\{1,2\}$, sends $e_1$ and $e_2$ to $ae_1$ and $e_2$ respectively for some $a\in{\bf k}^*$. It is easy to see that $a=1$ for $i=1$ and $a$ can be arbitrary for $i=2$.

5. Since ${\bf B}_3(V,V)$ is generated by $e_2$, any automorphism of ${\bf B}_3$ sends $e_1$ and $e_2$ to $ae_1+be_2$ and $ce_2$ respectively, where $a,c\in{\bf k}^*$ and $b\in{\bf k}$. It is easy to see that such a map is an automorphism iff $a=1$.

6. It follows from the proof of Lemma \ref{Cclass} that any automorphism of the algebra ${\bf C}(\alpha,\beta)$ sends $e_1$ and $e_2$ to $\pm e_1$ and $e_2$ respectively. It is easy to see that the map that sends $e_1$ and $e_2$ to $-e_1$ and $e_2$ respectively is an automorphism iff $\beta=0$ or $char{\bf k}=2$.

7. It follows from the proof of Lemma \ref{Dclass} that any automorphism of the algebra ${\bf D}_1(\alpha,\beta)$ sends $e_1$ to $e_1$ and sends $e_2$ either to $e_2$ or $e_1-e_2$. It is easy to see that the map that sends $e_1$ and $e_2$ to $e_1$ and $e_1-e_2$ respectively is an automorphism iff $\beta=2\alpha-1$.

8. It follows from the proof of Lemma \ref{Dclass} that any automorphism of the algebra ${\bf D}_2(\alpha,\beta)$ sends $e_1$ and $e_2$ to $e_1$ and $ae_2$ respectively for some $a\in{\bf k}^*$. It follows from the same proof that any automorphism of the algebra ${\bf D}_3(\alpha,\beta)$ is trivial.

9. It follows from the proof of Lemma \ref{Dclass} that any automorphism of ${\bf D}_4(\alpha,\beta)$ sends $e_1$ and $e_2$ to $e_1$ and $ae_2$, and $a=1$.
\end{proof}

We will finish the proof of Theorem \ref{alg} in the next section devoted to the algebras of the class ${\bf E}$.

\section{Algebras of the class ${\bf E}$}

In this section we consider the algebras of the class ${\bf E}$. It is clear that such an algebra is isomorphic to ${\bf E}_1(\Gamma)$ for some $\Gamma\in {\bf k}^4$.
Firstly, we describe isomorphisms inside this set and, thus, finish the proof of Theorem \ref{alg}.

\begin{lemma}\label{Eclass}
${\bf E}_1(\Gamma_1)\cong {\bf E}_1(\Gamma_2)$ iff one of the following conditions holds:
\begin{itemize}
\item $\Gamma_1=\Gamma_2$;
\item $\mathcal{C}_1(\Gamma_1)=\mathcal{C}_2(\Gamma_2)$ and $\mathcal{C}_2(\Gamma_1)=\mathcal{C}_1(\Gamma_2)$;
\item $\mathcal{C}_1(\Gamma_1), \mathcal{C}_1(\Gamma_2), \mathcal{C}_2(\Gamma_1), \mathcal{C}_2(\Gamma_2)\in\mathcal{T}$, $\mathcal{C}_1(\Gamma_1)\not=\mathcal{C}_2(\Gamma_1)$, $\mathcal{C}_1(\Gamma_2)\not=\mathcal{C}_2(\Gamma_2)$;
\item $\mathcal{C}_1(\Gamma_2),\mathcal{C}_2(\Gamma_2)\not\in \mathcal{T}$, $\mathcal{D}(\Gamma_2)\not=0$, and there is some $\sigma\in S_3$ such that $\mathcal{C}_i(\Gamma_1)=\mathcal{C}_{\sigma(i)}(\Gamma_2)$ for $i\in\{1,2,3\}$.
\end{itemize}
\end{lemma}
\begin{proof}
Suppose that $g\in GL(V)$ is such that $g*{\bf E}_1(\Gamma_1)={\bf E}_1(\Gamma_2)$. Then $ge_1$ and $ge_2$ are two linearly independent idempotents of ${\bf E}_1(\Gamma_2)$. Let us describe all nonzero idempotents of this algebra. Let $\Gamma_2=(\alpha,\beta,\gamma,\delta)$ and $u=xe_1+ye_2$ be some element of $V$. Then ${\bf E}_1(\Gamma_2)(u,u)=u$ iff
$x=x^2+(\alpha+\gamma)xy\mbox{ and }y=(\beta+\delta)xy+y^2.$
The solutions $(x,y)=(0,0),(1,0),(0,1)$ give the obvious  idempotents $0$, $e_1$ and $e_2$. All the other pairs $(x,y)$ satisfying the obtained equations are the solutions of the system of linear equations
$$
\begin{cases}
(\beta+\delta)x+y=1,\\
x+(\alpha+\gamma)y=1
\end{cases}
$$
with the additional conditions $x\not=0$ and $y\not=0$.
Let us consider the following cases:
\begin{itemize}
\item $\mathcal{C}_1(\Gamma_2),\mathcal{C}_2(\Gamma_2)\in \mathcal{T}$, i.e. $\alpha+\gamma=\beta+\delta=1$. In this case ${\bf E}_1(\Gamma_2)(u,u)=u$ iff either $u=0$ or $x+y=1$. Thus, $ge_1=ae_1+(1-a)e_2$ and $ge_2=be_1+(1-b)e_2$ for two different $a,b\in{\bf k}$. One can check that in this case
$$\Gamma_1=((1-b)\alpha+b\delta, a\beta+(1-a)\gamma,b\beta+(1-b)\gamma, (1-a)\alpha+a\delta).$$
If $\mathcal{C}_1(\Gamma_2)=\mathcal{C}_2(\Gamma_2)$, i.e. $(\beta,\delta)=(\gamma,\alpha)$, then we get $\Gamma_1=\Gamma_2$. If $\mathcal{C}_1(\Gamma_2)\not=\mathcal{C}_2(\Gamma_2)$, then the formula above gives all the possible $\Gamma_1$ with $\mathcal{C}_1(\Gamma_1),\mathcal{C}_2(\Gamma_1)\in \mathcal{T}$ and $\mathcal{C}_1(\Gamma_1)\not=\mathcal{C}_2(\Gamma_1)$.

\item one of the following three conditions holds:\\
1. $\mathcal{C}_1(\Gamma_2)\in \mathcal{T}$, $\mathcal{C}_2(\Gamma_2)\not\in \mathcal{T}$, i.e. $\beta+\delta=1$, $\alpha+\gamma\not=1$;\\
2. $\mathcal{C}_1(\Gamma_2)\not\in \mathcal{T}$, $\mathcal{C}_2(\Gamma_2)\in \mathcal{T}$, i.e. $\beta+\delta\not=1$, $\alpha+\gamma=1$;\\
3. $\mathcal{C}_1(\Gamma_2),\mathcal{C}_2(\Gamma_2)\not\in \mathcal{T}$, $\mathcal{D}(\Gamma_2)=0$, i.e. $\beta+\delta,\alpha+\gamma\not=1$ and our system has zero determinant.\\
It is easy to see that in all of these cases our system of linear equations doesn't have solutions satisfying the additional conditions, i.e. $e_1$ and $e_2$ are all the nonzero idempotents of ${\bf E}_1(\Gamma_2)$. Thus, either $ge_1=e_1$, $ge_2=e_2$ and $\Gamma_1=\Gamma_2$ or $ge_1=e_2$, $ge_2=e_1$, $\mathcal{C}_1(\Gamma_1)=\mathcal{C}_2(\Gamma_2)$ and $\mathcal{C}_2(\Gamma_1)=\mathcal{C}_1(\Gamma_2)$.

\item $\mathcal{C}_1(\Gamma_2),\mathcal{C}_2(\Gamma_2)\not\in \mathcal{T}$, $\mathcal{D}(\Gamma_2)\not=0$, i.e. $\beta+\delta,\alpha+\gamma\not=1$ and our system has nonzero determinant. In this case $(x,y)=\left(\frac{\alpha+\gamma-1}{\mathcal{D}(\Gamma_2)},\frac{\beta+\delta-1}{\mathcal{D}(\Gamma_2)}\right)$ is the unique solution of our system of linear equations. Hence, $e_1$, $e_2$ and $e_3=\frac{\alpha+\gamma-1}{\mathcal{D}(\Gamma_2)}e_1+\frac{\beta+\delta-1}{\mathcal{D}(\Gamma_2)}e_2$ are all the nonzero idempotents of ${\bf E}_1(\Gamma_2)$. Thus, there is $\sigma\in S_3$ such that $ge_i=e_{\sigma(i)}$ for $i\in\{1,2\}$. Then direct calculations show that $\mathcal{C}_i(\Gamma_1)=\mathcal{C}_{\sigma(i)}(\Gamma_2)$ for $i\in\{1,2,3\}$.
\end{itemize}
\end{proof}

Now we can finish the proof of Theorem \ref{alg}.

\begin{proof}[Proof of Theorem \ref{alg}] By Corollary \ref{classes}, the algebra $A$ belongs to one of the classes ${\bf A}$--${\bf E}$ and the class containing $A$ is unique. If $A$ belongs to one of the classes ${\bf A}$--${\bf D}$, then the statement of the theorem follows from Lemmas \ref{Aclass}--\ref{Dclass}.

Suppose that $A$ belongs to the class ${\bf E}$. Then $A$ can be represented by ${\bf E}_1(\Gamma)$ for some $\Gamma=(\alpha,\beta,\gamma,\delta)\in {\bf k}^4$. Now we have
\begin{itemize}
\item if $\mathcal{C}_1(\Gamma)=\mathcal{C}_2(\Gamma)\in \mathcal{T}$, then ${\bf E}_1(\Gamma)= {\bf E}_5(\beta)$;

\item if $\mathcal{C}_1(\Gamma),\mathcal{C}_2(\Gamma)\in \mathcal{T}$ and $\mathcal{C}_1(\Gamma)\not=\mathcal{C}_2(\Gamma)$, then ${\bf E}_1(\Gamma)\cong {\bf E}_1(1,1,0,0)= {\bf E}_4$ by Lemma \ref{Eclass};

\item if $\mathcal{C}_1(\Gamma)\in \mathcal{T}$ and $\mathcal{C}_2(\Gamma)\not\in \mathcal{T}$, then ${\bf E}_1(\Gamma)\cong {\bf E}_1(1-\beta,\gamma,\beta,\alpha)={\bf E}_2(\beta,\gamma,\alpha)$ by Lemma \ref{Eclass};

\item if $\mathcal{C}_1(\Gamma)\not\in \mathcal{T}$ and $\mathcal{C}_2(\Gamma)\in \mathcal{T}$, then ${\bf E}_1(\Gamma)={\bf E}_2(\gamma,\beta,\delta)$;

\item if $\mathcal{C}_1(\Gamma),\mathcal{C}_2(\Gamma)\not\in \mathcal{T}$, $\mathcal{D}(\Gamma)=0$ and $\alpha+\gamma\in{\bf k}^*_{>1}$, then ${\bf E}_1(\Gamma)={\bf E}_3\big(\gamma(\beta+\delta),\beta(\alpha+\gamma),\alpha+\gamma\big)$;

\item if $\mathcal{C}_1(\Gamma),\mathcal{C}_2(\Gamma)\not\in \mathcal{T}$, $\mathcal{D}(\Gamma)=0$ and $\alpha+\gamma\not\in{\bf k}^*_{>1}$, then ${\bf E}_1(\Gamma)\cong {\bf E}_1(\delta,\gamma,\beta,\alpha)={\bf E}_3\big(\beta(\alpha+\gamma),\gamma(\beta+\delta),\beta+\delta\big)$ by Lemma \ref{Eclass};

\item in the last two cases we have obtained a structure ${\bf E}_3(\alpha',\beta',\gamma')$ with $(\alpha',\beta',\gamma')\in{\bf k}^2\times{\bf k}^*_{>1}$. Interchanging $e_1$ and $e_2$, one gets ${\bf E}_3(\alpha',\beta',\gamma')\cong{\bf E}_3\big(\beta',\alpha',(\gamma')^{-1}\big)$; in particular, ${\bf E}_3(\alpha',\beta',-1)\cong{\bf E}_3(\beta',\alpha',-1)$. Thus, if $\gamma'=-1$ and $(\alpha',\beta')\not\in\mathcal{S}$, then we replace ${\bf E}_3(\alpha',\beta',-1)$ by ${\bf E}_3(\beta',\alpha',-1)$, and in all cases $A$ is represented by ${\bf E}_3(\Gamma')$ with $\Gamma'\in\mathcal{W}$;

\item if $\mathcal{C}_1(\Gamma),\mathcal{C}_2(\Gamma)\not\in \mathcal{T}$ and $\mathcal{D}(\Gamma)\not=0$, then there is a unique $\sigma\in S_3$ such that ${}^{\sigma^{-1}}(\mathcal{C}_1(\Gamma),\mathcal{C}_2(\Gamma),\mathcal{C}_3(\Gamma))\in\mathcal{\tilde V}$ and we have ${\bf E}_1(\Gamma)\cong {\bf E}_1(\Gamma')$ by Lemma \ref{Eclass}, where $\Gamma'\in\mathcal{V}$ is such that $\mathcal{C}_i(\Gamma')=\mathcal{C}_{\sigma(i)}(\Gamma)$ for $i\in\{1,2,3\}$.
\end{itemize}
By Lemma \ref{Eclass}, the structures from the set $$\{{\bf E}_1(\Gamma)\}_{\Gamma\in\mathcal{V}}\cup\{{\bf E}_2(\alpha,\beta,\gamma)\}_{(\alpha,\beta,\gamma)\in{\bf k}^3\setminus{\bf k}\times \mathcal{T}}\cup
\{{\bf E}_3(\alpha,\beta,\gamma)\}_{(\alpha,\beta,\gamma)\in\mathcal{W}}\cup\{{\bf E}_4\}\cup\{{\bf E}_5(\alpha)\}_{\alpha\in {\bf k}}$$ are pairwise non-isomorphic. Indeed, by the proof of Lemma \ref{Eclass}, for $(\alpha,\beta,\gamma),(\alpha',\beta',\gamma')\in{\bf k}^2\times{\bf k}^*_{>1}$ the algebras ${\bf E}_3(\alpha,\beta,\gamma)$ and ${\bf E}_3(\alpha',\beta',\gamma')$ are isomorphic iff either $(\alpha',\beta',\gamma')=(\alpha,\beta,\gamma)$ or $(\alpha',\beta',\gamma')=(\beta,\alpha,\gamma^{-1})$, and the second possibility with $\gamma,\gamma^{-1}\in{\bf k}^*_{>1}$ occurs only for $\gamma=-1$.
\end{proof}

As a consequence of the proof of Lemma \ref{Eclass} we can describe the automorphism groups of algebras of the class ${\bf E}$.

\begin{corollary}
$1.$ For $\Gamma\in\mathcal{V}$, ${\rm Aut}({\bf E}_1(\Gamma))$ is trivial if $\mathcal{C}_1(\Gamma)\not=\mathcal{C}_2(\Gamma)$, ${\rm Aut}({\bf E}_1(\Gamma))\cong C_2$ if $\mathcal{C}_1(\Gamma)=\mathcal{C}_2(\Gamma)\not=(-1,-1)$, ${\rm Aut}({\bf E}_1(-1,-1,-1,-1))\cong S_3$ if $char{\bf k}\not=3$.

$2.$ ${\rm Aut}({\bf E}_4)$ and ${\rm Aut}({\bf E}_2(\alpha,\beta,\gamma))$ are trivial for $(\alpha,\beta,\gamma)\in{\bf k}^3\setminus{\bf k}\times\mathcal{T}$.

$3.$ For $(\alpha,\beta,\gamma)\in {\bf k}^2\times{\bf k}^*_{>1}$, ${\rm Aut}({\bf E}_3(\alpha,\beta,\gamma))\cong C_2$ if $\gamma=-1$ and $\alpha=\beta$ and ${\rm Aut}({\bf E}_3(\alpha,\beta,\gamma))$ is trivial otherwise.

$4.$ ${\rm Aut}({\bf E}_5(\alpha))$ is isomorphic to the subgroup of $GL_2({\bf k})$ formed by matrices of the form $\left(\begin{array}{cc}
a&b\\
1-a&1-b
\end{array}\right)$, where $a,b\in{\bf k}$, $a\not=b$.

In particular, we have ${\rm dim}\,{\rm Aut}({\bf E}_1(\Gamma))={\rm dim}\,{\rm Aut}({\bf E}_2(\alpha,\beta,\gamma))={\rm dim}\,{\rm Aut}({\bf E}_3(\alpha,\beta,\gamma))={\rm dim}\,{\rm Aut}({\bf E}_4)=0$ and ${\rm dim}\,{\rm Aut}({\bf E}_5(\alpha))=2$.
\end{corollary}
\begin{proof}
1. Any automorphism of ${\bf E}_1(\Gamma)$ has to send $e_1$ and $e_2$ to $e_{\sigma(1)}$ and $e_{\sigma(2)}$ respectively for some $\sigma\in S_3$, where $e_3$ is defined in the proof of Lemma \ref{Eclass}. Such a map is an automorphism iff $\mathcal{C}_i(\Gamma)=\mathcal{C}_{\sigma(i)}(\Gamma)$ for $i=1,2$. If $\mathcal{C}_1(\Gamma)\not=\mathcal{C}_2(\Gamma)$, then we have also $\mathcal{C}_3(\Gamma)\not=\mathcal{C}_1(\Gamma),\mathcal{C}_2(\Gamma)$ and, hence, only identical element of $S_3$ determines an automorphism. If $\mathcal{C}_1(\Gamma)=\mathcal{C}_2(\Gamma)\not=(-1,-1)$, then one can check that $\mathcal{C}_3(\Gamma)\not=\mathcal{C}_1(\Gamma),\mathcal{C}_2(\Gamma)$ and, hence, only identical element of $S_3$ and the element that swaps $e_1$ and $e_2$ determine automorphisms. If $\mathcal{C}_1(\Gamma)=\mathcal{C}_2(\Gamma)=(-1,-1)$, then any $\sigma\in S_3$ determines an automorphism. Note that $(-1,-1,-1,-1)\in\mathcal{V}$ iff $char{\bf k}\not=3$.

2. Follows directly from the proof of Lemma \ref{Eclass}.

3. It follows from the proof of Lemma \ref{Eclass} that an automorphism of ${\bf E}_3(\alpha,\beta,\gamma)$ is either  trivial or swaps $e_1$ and $e_2$. The last mentioned map is an automorphism iff $\gamma=-1$ and $\alpha=\beta$.

4. It follows directly from the proof of Lemma \ref{Eclass} that automorphisms of ${\bf E}_5(\alpha)$ are exactly the linear maps that send $e_1$ and $e_2$ to $ae_1+(1-a)e_2$ and $be_1+(1-b)e_2$ for two different $a,b\in{\bf k}$.
\end{proof}

Now we are going to discuss some facts about degenerations of the form $A\rightarrow B$, where $A$ is an algebra of the class ${\bf E}$. First of all, let us prove the following lemma.

\begin{lemma}\label{EtoD}
$1.$ For any $\Gamma\in \mathcal{V}$ and $(\beta,\gamma)\in\mathcal{C}(\Gamma)$ there exists a degeneration ${\bf E}_1(\Gamma)\to {\bf D}_2(\beta,\gamma)$.\\
$2.$ For any $(\alpha,\beta,\gamma)\in {\bf k}^3\setminus {\bf k}\times\mathcal{T}$ there exists a degeneration ${\bf E}_2(\alpha,\beta,\gamma)\to {\bf D}_2(\beta,\gamma)$.\\
$3.$ For any $(\alpha,\delta,\epsilon)\in {\bf k}^2\times{\bf k}^*_{>1}$ and $(\beta,\gamma)\in\mathcal{C}(\alpha,\delta,\epsilon)$ there exists a degeneration ${\bf E}_3(\alpha,\delta,\epsilon)\to {\bf D}_2(\beta,\gamma)$.
\end{lemma}
\begin{proof}
The parametrized basis $E_1^t=e_1$, $E_2^t=te_2$ gives the degeneration ${\bf E}_1(\Gamma)\to {\bf D}_2\big(\mathcal{C}_1(\Gamma)\big)$ for any $\Gamma\in {\bf k}^4$. If $\Gamma\in\mathcal{V}$ and $1\le i\le 3$, then, by Lemma \ref{Eclass}, there exists $\Gamma'\in{\bf k}^4$ such that ${\bf E}_1(\Gamma)\cong {\bf E}_1(\Gamma')$ and $\mathcal{C}_i(\Gamma)=\mathcal{C}_1(\Gamma')$. Hence, ${\bf E}_1(\Gamma)\cong {\bf E}_1(\Gamma')\to {\bf D}_2\big(\mathcal{C}_i(\Gamma)\big)$. We also have
${\bf E}_2(\alpha,\beta,\gamma)={\bf E}_1(1-\alpha,\beta,\alpha,\gamma)\to {\bf D}_2(\beta,\gamma)$ for $(\alpha,\beta,\gamma)\in {\bf k}^3\setminus {\bf k}\times\mathcal{T}$ and
{\small\begin{multline*}{\bf E}_3(\alpha,\delta,\epsilon)={\bf E}_1\left((1-\alpha)\epsilon,\frac{\delta}{\epsilon},\alpha\epsilon,\frac{1-\delta}{\epsilon}\right)\cong{\bf E}_1\left(\frac{1-\delta}{\epsilon},\alpha\epsilon,\frac{\delta}{\epsilon},(1-\alpha)\epsilon\right)
\to
{\bf D}_2\big(\alpha\epsilon,(1-\alpha)\epsilon\big),{\bf D}_2\left(\frac{\delta}{\epsilon},\frac{1-\delta}{\epsilon}\right)\end{multline*}}
for $(\alpha,\delta,\epsilon)\in {\bf k}^2\times{\bf k}^*_{>1}$.
\end{proof}

Let us define, for $\Gamma=(\alpha,\beta,\gamma,\delta)\in{\bf k}^4$, the following subset of $\mathcal{A}_2$:
{\footnotesize $$\mathcal{G}(\Gamma)=
\left\{ \mu\left|  \begin{array} {@{}L{0.79\textwidth}@{}} 	c_{22}^1=0;\,\,c_{21}^1=\gamma c_{22}^2;\,\,c_{12}^1=\alpha c_{22}^2;\\
\big(1-\gamma-\delta(\alpha+\gamma)\big)c_{12}^2-\big(1-\alpha-\beta(\alpha+\gamma)\big)c_{21}^2=  \ \big(\beta(1-\gamma)-\delta(1-\alpha)\big)c_{11}^1;\\
  \big(1-\alpha-\beta(\alpha+\gamma)\big)^2c_{11}^2c_{22}^2=  (\beta c_{11}^1-c_{12}^2)\big(\mathcal{D}(\Gamma)c_{12}^2+((\alpha-1)(\delta-1)-\beta\gamma)c_{11}^1\big);\\
 \big(1-\gamma-\delta(\alpha+\gamma)\big)^2c_{11}^2c_{22}^2=  (\delta c_{11}^1-c_{21}^2)\big(\mathcal{D}(\Gamma)c_{21}^2+((\beta-1)(\gamma-1)-\alpha\delta)c_{11}^1\big);\\
 \big(1-\alpha-\beta(\alpha+\gamma)\big)\big(1-\gamma-\delta(\alpha+\gamma)\big)c_{11}^2c_{22}^2=  \  (\beta c_{11}^1-c_{12}^2)\big(\mathcal{D}(\Gamma)c_{21}^2+((\beta-1)(\gamma-1)-\alpha\delta)c_{11}^1\big);\\
 \big(1-\alpha-\beta(\alpha+\gamma)\big)\big(1-\gamma-\delta(\alpha+\gamma)\big)c_{11}^2c_{22}^2= \ (\delta c_{11}^1-c_{21}^2)\big(\mathcal{D}(\Gamma)c_{12}^2+((\alpha-1)(\delta-1)-\beta\gamma)c_{11}^1\big)\\
\end{array} \right.\right\}. $$}
Here and further in a definition of a subset of $\mathcal{A}_2$ we always assume by default that $\mu\in\mathcal{A}_2$ and $c_{ij}^k$ ($i,j,k\in\{1,2\})$ are structure constants of $\mu$.
The following lemma will allow us to use $\mathcal{G}(\Gamma)$ as a separating set for some non-degenerations. Its proof is a direct calculation and so it is left to the reader.

\begin{lemma}
The set $\mathcal{G}(\Gamma)$ is closed upper invariant and contains ${\bf E}_1(\Gamma)$ for any $\Gamma\in{\bf k}^4$.
\end{lemma}

\section{Degenerations of $2$-dimensional algebras}

In this section we describe all degenerations of $2$-dimensional algebras. Note that the results are valid for algebras over an algebraically closed field of arbitrary characteristic.






\begin{theorem} \label{graf2}
$\mathcal{A}_2$ has the graph of primary degenerations presented in Figure 1.

\end{theorem}

\begin{proof} All primary degenerations that don't follow from Lemma \ref{EtoD} are presented in Table 2.
Table 3 describes separating sets for all required non-degenerations and, thus, finishes the proof of the theorem.

The verification of degenerations is an easy direct calculation in each case. An example clarifying how to do this can be found in the proof of \cite[Theorem 3]{kppv}.
The verification of Table 3 is more difficult. To clarify how one can fulfill it, let us consider the row ${\bf A}_4(\alpha)\not\to{\bf B}_2(\gamma),{\bf D}_2(\beta,\gamma),{\bf E}_5(\beta)$ of the table. It is easy to prove that ${\bf A}_4(\alpha)$ belongs to the presented separating set that we denote by $\mathcal{R}$.
The fact that $\mathcal{R}$ is upper invariant can be checked by a direct calculation. What exactly one has to check is explained in Section 2.
Let us prove that the orbits of ${\bf B}_2(\gamma)$, ${\bf D}_2(\beta,\gamma)$, and ${\bf E}_5(\beta)$ don't intersect $\mathcal{R}$. Let $\lambda$ be one of these structures. Suppose that the structure constants $c_{ij}^k$ ($i,j,k=1,2$) of $\lambda$ in the basis $f_1$, $f_2$ satisfy the defining equations of $\mathcal{R}$. Then $c_{22}^1=c_{22}^2=0$, and hence $\lambda(f_2,f_2)=0$. Since $\mathcal{R}$ is invariant under the basis rescaling, we may assume that $f_2\in\{e_1,e_2\}$ if $\lambda={\bf B}_2(\gamma)$, $f_2=e_2$ if $\lambda={\bf D}_2(\beta,\gamma)$, and $f_2=e_2-e_1$ if $\lambda={\bf E}_5(\beta)$. Now, in view of the  upper invariance of $\mathcal{R}$, we may assume that $(f_1,f_2)\in\{(e_1,e_2),(e_2,e_1)\}$ if $\lambda={\bf B}_2(\gamma)$, $(f_1,f_2)=(e_1,e_2)$ if $\lambda={\bf D}_2(\beta,\gamma)$, and $(f_1,f_2)=(e_1,e_2-e_1)$ if $\lambda={\bf E}_5(\beta)$.
We have 
\begin{itemize}
    \item $c_{12}^1+c_{21}^1=1\not=0$ if $\lambda={\bf B}_2(\gamma)$, $f_1=e_1$, $f_2=e_2$;
    \item $c_{12}^2+c_{21}^2=1\not=0=c_{11}^1$ if $\lambda={\bf B}_2(\gamma)$, $f_1=e_2$, $f_2=e_1$;
    \item $\alpha^2c_{12}^1c_{11}^2=0\not=(c_{11}^1)^2$ if either $\lambda={\bf D}_2(\beta,\gamma)$, $f_1=e_1$, $f_2=e_2$ or $\lambda={\bf E}_5(\beta)$, $f_1=e_1$, $f_2=e_2-e_1$.
\end{itemize}
 Thus, the structure constants of $\lambda$ in any basis don't satisfy the defining equations of $\mathcal{R}$, i.e. $O(\lambda)\cap\mathcal{R}=\varnothing$. The other nondegenerations can be considered in the same manner. Note that the structure ${\bf B}_2(\gamma)$ does not lie in the separating set given for it in Table 3, but the structure obtained from ${\bf B}_2(\gamma)$ by interchanging $e_1$ and $e_2$ does.
\end{proof}

Let us recall that $n$-dimensional algebra $A$ has a {\it level} $m$ if
\begin{itemize}
\item there exists a sequence of $n$-dimensional algebras $A_0,\dots, A_m$ such that $A_0={\bf k}^n$, $A_m=A$ and, for $0\le i\le m-1$, one has $A_{i+1}\to A_i$ and $A_{i+1}\not\cong A_i$;
\item if $A_0,\dots, A_{m+1}$ is a sequence of algebras such that $A_0={\bf k}^n$, $A_{m+1}=A$, and $A_{i+1}\to A_i$ for $0\le i\le m$, then $A_{i+1}\cong A_i$ for some $0\le i\le m$.
\end{itemize}
Theorem \ref{graf2} gives us the following partition of $\mathcal{A}_2$ to levels:
\begin{center}
\small\setlength{\tabcolsep}{4pt}
\begin{longtable}{|L{0.13\textwidth}|L{0.05\textwidth}|L{0.14\textwidth}|L{0.21\textwidth}|L{0.30\textwidth}|}
\hline
\mbox{level}  &  0 &  1 &  2 &  3 \\
\hline
\begin{array}{l}
\mbox{algebra} \\
\mbox{structures}
\end{array} & {\bf k}^2 &
\begin{array}{l}{\bf A}_3, {\bf B}_3,  {\bf E}_5(\alpha) \end{array} &
\begin{array}{l}
{\bf A}_1(\alpha), {\bf A}_2,   {\bf B}_2(\alpha),\\ {\bf D}_2(\alpha,\beta),  {\bf E}_4
\end{array}&\begin{array}{l}{\bf A}_4(\alpha), {\bf B}_1(\alpha), {\bf C}(\alpha,\beta), \\
{\bf D}_1(\alpha,\beta), {\bf D}_3(\alpha,\beta),\\ {\bf D}_4(\alpha,\beta), {\bf E}_1(\Gamma),\\  {\bf E}_2(\alpha,\beta,\gamma), {\bf E}_3(\alpha,\beta,\gamma)\end{array}\\
\hline
\end{longtable}
\end{center}
Note also that the algebras ${\bf k}^2$, ${\bf B}_3$, ${\bf E}_4$ and ${\bf E}_5(\alpha)$ ($\alpha\in{\bf k}$) form a closed subset of $\mathcal{A}_2$ that has two interesting descriptions. First of all, these are exactly all $2$-dimensional algebras that don't degenerate to ${\bf A}_3$. Secondly, any one-generated subalgebra of such an algebra is $1$-dimensional and this property doesn't hold for other algebras. In fact, these two descriptions define the same set of algebras in a variety of algebras of any dimension. Note also that ${\bf E}_4$ is a unique $2$-dimensional algebra of level two that doesn't have non-trivial derivations.

\section{Closures for orbits of infinite series}

In this section we describe closures of orbits for infinite series from our classification. To make this description nicer and more complete we introduce two additional series and one additional algebra. For $\alpha\in{\bf k}$, we introduce ${\bf D}_2'(\alpha)={\bf D}_2(\alpha,-\alpha)$ and ${\bf D}_3'(\alpha)={\bf D}_3(\alpha,-\alpha)$. Note that ${\bf D}_2'(*)\subset {\bf D}_2(*)$ and ${\bf D}_3'(*)\subset {\bf D}_3(*)$. Also we define ${\bf A}_4'={\bf A}_4(0)\in {\bf A}_4(*)$. Here and further, for a symbol $\bf X$, we denote by ${\bf X}(*)$ the set formed by all ${\bf X}(\Gamma)$ that are defined. For example, ${\bf D}_2(*)=\{{\bf D}_2(\Gamma)\mid \Gamma\in{\bf k}^2\}$, ${\bf E}_3(*)=\{{\bf E}_3(\Gamma)\mid \Gamma\in {\bf k}^2\times{\bf k}^*\}$.

\begin{theorem}\label{inf2} For each row of Table 4, the second column contains all isomorphism classes of algebras whose orbits lie in the closure of the orbit of the series of algebras contained in the first column of the same row.
\end{theorem}
\begin{proof} 
All required degenerations that don't follow from Theorem \ref{graf2} are proved in Table 5.
Table 6 describes separating sets for all required non-degenerations and, thus, finishes the proof of the theorem.

Note that the structures ${\bf B}_1(\alpha)$, ${\bf D}_1(\alpha,\beta)$, and ${\bf E}_3(\alpha,\beta,\gamma)$ don't lie in the corresponding separating sets, but the structures
$$
\left(\begin{array}{cc}
0     &1  \\
1     &0 
\end{array}\right)*{\bf B}_1(\alpha),\,\,
\left(\begin{array}{cc}
0     &1  \\
1     &0 
\end{array}\right)*{\bf D}_1(\alpha,\beta),\,\,\mbox{ and }
\left(\begin{array}{cc}
1     &\gamma  \\
0     &1 
\end{array}\right)*{\bf E}_3(\alpha,\beta,\gamma)
$$
satisfy the required conditions.

In Table 5 the degeneration ${\bf E}_3(*)\to {\bf D}_3(\alpha,\beta)$ is given only for $\alpha+\beta\neq 0$ and the degeneration ${\bf D}_3(*)\to {\bf A}_4(\alpha)$ only for $\alpha\neq 0$ (for $\alpha+\beta=0$, resp. $\alpha=0$, the parametrized bases and indices given there are not defined). The remaining cases follow from the other rows of Table 5: since ${\bf E}_3(*)\to {\bf D}_1(\alpha,\beta)$ for all $(\alpha,\beta)$, we have $\overline{\bigcup_{\alpha,\beta}O({\bf D}_1(\alpha,\beta))}\subset\overline{\bigcup_{\Gamma}O({\bf E}_3(\Gamma))}$, and hence ${\bf D}_1(*)\to {\bf D}_3'(\alpha)$ implies ${\bf E}_3(*)\to {\bf D}_3(\alpha,-\alpha)$; analogously, ${\bf D}_3'(*)\to {\bf A}_4'$ and ${\bf D}_3'(*)\subset{\bf D}_3(*)$ imply ${\bf D}_3(*)\to {\bf A}_4(0)$.

The structure ${\bf B}_1(\alpha)$ does not lie in the set $\{c_{22}^1=c_{22}^2=c_{12}^1+c_{21}^1=0\}$ of Table~6, but the structure obtained from it by interchanging $e_1$ and $e_2$ does. The orbit of ${\bf D}_4(\alpha,\beta)$ does not meet this set, nor the sets of Table~6 for ${\bf C}(*)$ and ${\bf D}_1(*)$: the only square zero elements of ${\bf D}_4(\alpha,\beta)$ are the multiples of $e_2$, and $e_1e_2+e_2e_1=e_1+e_2$.

The row ${\bf D}_2(*)\not\to{\bf A}_4(\alpha),{\bf B}_1(\alpha),{\bf D}_3(\alpha,\beta),{\bf E}_4$ of Table~6 is needed because all the other separating sets of Table~6 containing ${\bf D}_2(*)$ contain also ${\bf A}_4(*)$. The set $\{c_{22}^1=c_{22}^2=c_{12}^1=c_{21}^1=0\}$ consists of the structures such that $e_2e_2=0$ and $\langle e_2\rangle$ is a two-sided ideal; it is closed upper invariant and contains ${\bf D}_2(\alpha,\beta)$, while the algebras ${\bf A}_4(\alpha)$, ${\bf B}_1(\alpha)$, ${\bf D}_3(\alpha,\beta)$ and ${\bf E}_4$ do not have a nonzero element $x$ such that $x^2=0$ and $xA+Ax\subset\langle x\rangle$. The remaining non-degenerations follow from Theorem \ref{graf2}, Table~6 and the already proved inclusions between the closures (for example, ${\bf E}_5(*)\not\to{\bf A}_3$ because ${\bf E}_4\to{\bf E}_5(\alpha)$ for all $\alpha$ and ${\bf E}_4\not\to{\bf A}_3$; ${\bf A}_4(*)\not\to{\bf C}(\alpha,\beta)$ because ${\bf D}_4(*)\to{\bf A}_4(\alpha)$ for all $\alpha$ and ${\bf D}_4(*)\not\to{\bf C}(\alpha,\beta)$).
\end{proof}

Theorems \ref{graf2} and \ref{inf2} give a lattice of subsets for $\mathcal{A}_2$. This lattice is presented in Figure 2.
In this figure the leftmost set coincides with $\mathcal{A}_2$ and sets placed in one column have the same dimension equal to the number standing above them.
Two sets of dimensions $i$ and $i+1$ are connected by an edge iff the set of dimension $i+1$ contains the set of dimension $i$. Moreover, if $X,Y\subset \mathcal{A}_2$ correspond to two vertices of the diagram, then $X\cap Y$ is equal to the union of all $Z\subset \mathcal{A}_2$ corresponding to vertices of the diagram such that there exist paths from $X$ to $Z$ and from $Y$ to $Z$ going from left to right. For example,
$\overline{O\big({\bf D}_1(*)\big)}\cap \overline{O\big({\bf C}(*)\big)}=\overline{O\big({\bf A}_4'\big)}$ and $\overline{O\big({\bf D}_1(*)\big)}\cap \overline{O\big({\bf D}_2(*)\big)}=\overline{O\big({\bf B}_2(*)\big)}\cup\overline{O\big({\bf D}_2'(*)\big)}$.














\section{Subvarieties Defined by Identities}

Now we are going to apply the results of previous sections to develop the varieties of two-dimensional flexible and bicommutative algebras. In particular, we will describe the varieties of commutative and anticommutative algebras. Since there exists only one nontrivial two-dimensional anticommutative algebra, the last mentioned problem is not of big interest. Note also that in the same way one can recover the results of \cite{BB09}, where the analogous problems were solved for two-dimensional Novikov and pre-Lie algebras. Since the classifications of flexible and bicommutative algebras depend on the characteristic of the ground field, we assume everywhere in this section that $char{\bf k}\not=2$.

\subsection{Flexible Algebras} By definition, an algebra is called \textit{flexible} if it satisfies the identity $(xy)x=x(yx).$ It is clear that all commutative and anticommutative algebras are flexible.
Using Theorem \ref{alg}, one can verify that any two-dimensional flexible algebra is either (anti)commutative or ${\bf E}_5(\alpha).$
For $\alpha,\beta\in{\bf k}$, let us introduce the algebras 
\begin{equation*}
\begin{aligned}{\bf D}_2^c(\alpha)&={\bf D}_2(\alpha,\alpha),&\quad{\bf E}_2^c(\alpha)&={\bf E}_2(\tfrac{1}{2},  \alpha, \alpha), \\
{\bf E}_3^c(\alpha)&={\bf E}_3(\tfrac{1}{2}, \tfrac{1}{2}, \alpha), &\quad 
{\bf E}_1^c(\alpha, \beta)&={\bf E}_1(\alpha, \beta, \alpha, \beta).
\end{aligned}
\end{equation*}
It follows from our classification that any nontrivial two-dimensional commutative algebra can be represented by a unique structure from the set
\begin{multline*}
\bigl\{{\bf A}_1(\tfrac{1}{2}),{\bf A}_3,{\bf B}_2(\tfrac{1}{2}),{\bf C}(\tfrac{1}{2},0),{\bf D}_1( \tfrac{1}{2},0),{\bf D}_4(\tfrac{1}{2},\tfrac{1}{2}),{\bf E}_5( \tfrac{1}{2})\bigr\}
\\
\cup \{{\bf D}_2^c(\alpha),{\bf E}_2^c(\alpha)\}_{\alpha\in{\bf k}\setminus\{\frac{1}{2}\} }\cup 
\{{\bf E}_3^c(\alpha)\}_{\alpha\in{\bf k}^*_{>1}}\cup \{{\bf E}_1^c(\alpha, \beta)\}_{(\alpha,\beta,\alpha,\beta) \in  \mathcal{V}}
\end{multline*}
and that any noncommutative flexible algebra
can be represented by a structure from the set $\{{\bf B}_3\}\cup\{{\bf E}_5(\alpha)\}_{\alpha\in{\bf k}\setminus\{\frac{1}{2}\} }$.
Using Theorem \ref{graf2}, we get the graph of primary degenerations for the variety of two-dimensional flexible algebras presented in Figure 3.

It is easy to see that the variety of commutative algebras is simply ${\bf k}^6$, {\it i.e.}, irreducible. 
Then it is clear that the variety of flexible algebras has two irreducible components.
The first component is  
$\overline{\{O({\bf E}_5(\alpha)) \}_{\alpha \in \bf k}}= 
\{ {\bf E}_5(\alpha), {\bf B}_3, {\bf k}^2 \}_{\alpha \in {\bf k}}$.
The second component, formed by all commutative algebras, is equal to the closure of the orbit of the algebra series ${\bf E}_1^c(*)$. 

This
variety of flexible algebras  does not have rigid algebras, and has the lattice of subsets presented in Figure 4.
The lattice satisfies the same properties as the lattice from the previous section. To prove this it is enough to use Theorem \ref{inf2} and its proof.
The only difference is that one must use the parametrized indices
$(\frac{1}{2}+\frac{t}{2}, \frac{1}{2}+\frac{t}{2} )$, $(\frac{1}{2}, \frac{1}{2}- \frac{t^2}{2}, \frac{1}{2}-\frac{t^2}{2} )$, and $(\frac{1}{2}, \frac{1}{2}, \frac{1}{t} )$ in the degenerations ${\bf D}_2(*) \to {\bf A}_1(\tfrac{1}{2})$, ${\bf E}_2(*) \to {\bf C}(\tfrac{1}{2},0)$, and ${\bf E}_3(*) \to {\bf D}_1(\tfrac{1}{2},0)$, respectively, to obtain the degenerations ${\bf D}_2^c(*) \to {\bf A}_1(\tfrac{1}{2})$, ${\bf E}_2^c(*) \to {\bf C}(\tfrac{1}{2},0)$, and ${\bf E}_3^c(*) \to {\bf D}_1(\tfrac{1}{2},0)$. Finally, the index $\big(\tfrac{1}{2},\tfrac{1}{2t},\tfrac{1}{2t}\big)$ in the degeneration ${\bf E}_2(*)\to{\bf D}_4(\alpha,\beta)$ of Table~5 gives ${\bf E}_2^c(*)\to{\bf D}_4(\tfrac{1}{2},\tfrac{1}{2})$, and the index $\big(\tfrac{1}{2},\tfrac{1}{2},\tfrac{1}{1-t}\big)$ with the basis $E_1^t=e_1$, $E_2^t=\frac{e_1-(1-t)e_2}{t}$ gives ${\bf E}_3^c(*)\to{\bf D}_4(\tfrac{1}{2},\tfrac{1}{2})$; the non-degenerations of ${\bf D}_4(\tfrac{1}{2},\tfrac{1}{2})$ follow from Table~3.

\subsection{Bicommutative Algebras}
The variety of bicommutative algebras (see, for example, \cite{Dzhuma11}) is defined by the identities
$x(yz)=y(xz) \mbox{ and } (xy)z=(xz)y.$ It follows from Theorem \ref{alg} that any nontrivial bicommutative algebra is isomorphic to a unique algebra from the set
$$
\{{\bf A}_3,{\bf B}_2(0),{\bf B}_2(1),{\bf D}_1(0,0),{\bf D}_2(1,1),{\bf D}_2(0,0),{\bf E}_1(0,0,0,0)\}.
$$
Using Theorem \ref{graf2}, we get the graph of primary degenerations  for the variety of two-dimensional bicommutative algebras. This graph is presented in Figure 5.

Thus, the irreducible components in the variety of two-dimensional bicommutative algebras are
\begin{align*}
\overline{O\big({\bf D}_1(0,0)\big)}&=\{{\bf D}_1(0,0),{\bf D}_2(0,0),{\bf B}_2(0),{\bf B}_2(1),{\bf A}_3,{\bf k}^2\},\\
\overline{O\big({\bf E}_1(0,0,0,0)\big)}&=\{{\bf E}_1(0,0,0,0),{\bf D}_2(0,0),{\bf D}_2(1,1),{\bf A}_3,{\bf k}^2\}.
\end{align*}
These components are generated by the rigid bicommutative algebras ${\bf D}_1(0,0)$ and ${\bf E}_1(0,0,0,0)$ and all have dimension $4$.

\section*{Appendix A. Tables.}
\addcontentsline{toc}{section}{Appendix A. Tables}
\begingroup\footnotesize\setlength{\tabcolsep}{4pt}

\begin{longtable}{|>{\raggedright\arraybackslash}p{0.24\textwidth}|>{\raggedright\arraybackslash}p{0.14\textwidth}|>{\raggedright\arraybackslash}p{0.16\textwidth}|>{\raggedright\arraybackslash}p{0.16\textwidth}|>{\raggedright\arraybackslash}p{0.14\textwidth}|}
\hline
\multicolumn{5}{|l|}{\textbf{Table 1}}\\
\hline
\endfirsthead
\hline
\multicolumn{5}{|l|}{\textbf{Table 1} (continued)}\\
\hline
\endhead
\hline\endfoot
\endlastfoot
\hline
${\bf A}_1(\alpha)$, $\alpha\in{\bf k}$ & $e_1e_1=e_1+e_2,$& $e_1e_2=\alpha e_2,$& $e_2e_1= (1-\alpha) e_2,$& $e_2e_2=0$ \\
\hline
${\bf A}_2$ & $e_1e_1=e_2,$& $e_1e_2=e_2,$& $e_2e_1= -e_2,$& $e_2e_2=0$ \\
\hline
${\bf A}_3$ & $e_1e_1=e_2,$ & $e_1e_2= 0,$ & $e_2e_1=0,$ & $e_2e_2=0$ \\
\hline
${\bf A}_4(\alpha)$, $\alpha\in{\bf k_{\ge 0}}$ & $e_1e_1= \alpha e_1+  e_2,$ & $e_1e_2= e_1+ \alpha e_2,$ & $ e_2e_1= - e_1,$ & $ e_2e_2=0$\\
\hline
${\bf B}_1(\alpha)$, $\alpha\in{\bf k}$ & $e_1e_1=0,$ & $e_1e_2=(1-\alpha)e_1 +  e_2,$ & $ e_2e_1=\alpha  e_1 - e_2,$ & $ e_2e_2=0$\\
\hline
${\bf B}_2(\alpha)$, $\alpha\in{\bf k}$ & $e_1e_1=0,$ & $e_1e_2=(1-\alpha)e_1,$ & $ e_2e_1=\alpha  e_1,$ & $ e_2e_2=0$\\
\hline
${\bf B}_3$ & $e_1e_1=0,$ & $ e_1e_2= e_2,$ & $e_2e_1=-e_2,$ & $ e_2e_2=0$ \\
\hline
${\bf C}(\alpha,\beta)$, $(\alpha,\beta)\in{\bf k}\times{\bf k_{\ge 0}}$ & $e_1e_1=e_2,$ & $ e_1e_2= (1-\alpha)e_1+\beta e_2,$ & $e_2e_1=\alpha e_1-\beta e_2,$ & $ e_2e_2=e_2$ \\
\hline
${\bf D}_1(\alpha,\beta)$, $(\alpha,\beta)\in \mathcal{U}$ & $e_1e_1=e_1,$ & $ e_1e_2=(1-\alpha)e_1+ \beta  e_2,$ & $e_2e_1=\alpha  e_1 - \beta   e_2,$ & $ e_2e_2=0$\\
\hline
${\bf D}_{2}(\alpha,\beta)$, $(\alpha,\beta)\in{\bf k}^2\setminus \mathcal{T}$ & $e_1e_1=e_1,$ & $ e_1e_2=\alpha e_2,$ & $e_2e_1=\beta e_2,$ & $ e_2e_2=0$\\
\hline
${\bf D}_{3}(\alpha,\beta)$, $(\alpha,\beta)\in{\bf k}^2\setminus \mathcal{T}$ & $e_1e_1=e_1,$ & $ e_1e_2=e_1+ \alpha  e_2,$ & $e_2e_1=- e_1 + \beta  e_2,$ & $ e_2e_2=0$\\
\hline
${\bf D}_{4}(\alpha,\beta)$, $(\alpha,\beta)\in{\bf k}^2$ & $e_1e_1=e_1,$ & $ e_1e_2=(1-\alpha)e_1+ \beta  e_2,$ & $e_2e_1=\alpha e_1 + (1-\beta)  e_2,$ & $ e_2e_2=0$\\
\hline
${\bf E}_1(\alpha,\beta,\gamma,\delta)$, $(\alpha,\beta,\gamma,\delta)\in \mathcal{V}$ & $e_1e_1=e_1,$ & $ e_1e_2=\alpha e_1+ \beta  e_2,$ & $e_2e_1=\gamma  e_1 + \delta  e_2,$ & $ e_2e_2=e_2$\\
\hline
$\begin{array}{l}
\!\!\!{\bf E}_2(\alpha,\beta,\gamma),\\
(\alpha,\beta,\gamma)\in {\bf k}^3\setminus {\bf k}\times\mathcal{T}\end{array}$ & $e_1e_1=e_1,$ & $ e_1e_2=(1-\alpha) e_1+ \beta  e_2,$ & $e_2e_1=\alpha  e_1 + \gamma  e_2,$ & $ e_2e_2=e_2$\\
\hline
$\begin{array}{l}
\!\!\!{\bf E}_3(\alpha,\beta,\gamma),\\(\alpha,\beta,\gamma)\in \mathcal{W}
\end{array}$ & $e_1e_1=e_1,$ & $ e_1e_2=(1-\alpha)\gamma e_1+ \frac{\beta}{\gamma}  e_2,$ & $e_2e_1=\alpha\gamma  e_1 + \frac{1-\beta}{\gamma}  e_2,$ & $ e_2e_2=e_2$\\
\hline
${\bf E}_4$& $e_1e_1=e_1,$ & $ e_1e_2=e_1+ e_2,$ & $e_2e_1=0,$ & $ e_2e_2=e_2$\\
\hline
${\bf E}_5(\alpha)$, $\alpha\in{\bf k}$& $e_1e_1=e_1,$ & $ e_1e_2=(1-\alpha) e_1+ \alpha  e_2,$ & $e_2e_1=\alpha  e_1 + (1-\alpha)  e_2,$ & $ e_2e_2=e_2$\\
\hline
\end{longtable}

\begin{longtable}{|L{0.36\textwidth}|L{0.27\textwidth}|L{0.27\textwidth}|}
\hline
\multicolumn{3}{|l|}{\textbf{Table 2}}\\
\hline
\mbox{degenerations} & \mbox{parametrized bases} & \\
\hline\hline
\endfirsthead
\hline
\multicolumn{3}{|l|}{\textbf{Table 2} (continued)}\\
\hline
\mbox{degenerations} & \mbox{parametrized bases} & \\
\hline\hline
\endhead
\hline\endfoot
\endlastfoot
{\bf A}_1(\alpha) \to  {\bf A}_3   &  E_1^t=te_1 & E_2^t=t^2e_2 \\
\hline
{\bf A}_1(\alpha) \to  {\bf E}_5(\alpha)   &  E_1^t=e_1 & E_2^t=e_1+t^{-1}e_2 \\
\hline
{\bf A}_2  \to  {\bf A}_3 &   E_1^t=te_1 & E_2^t=t^2e_2 \\
\hline
{\bf A}_2  \to  {\bf B}_3 &   E_1^t=e_1 & E_2^t= t^{-1}e_2\\
\hline
{\bf A}_4(\alpha)  \to  {\bf A}_2   &   E_1^t=te_1-e_2 & E_2^t= t^2e_2 \\
\hline
{\bf B}_1(\gamma)  \to  {\bf A}_2   &  
E_1^t=e_1+te_2 & E_2^t=-t^2e_2 \\
\hline
{\bf B}_1(\gamma)  \to  {\bf B}_2(\gamma)  &  
E_1^t=te_1 & E_2^t=e_2 \\
\hline
{\bf B}_2(\gamma)  \to  {\bf A}_3 &  
E_1^t=e_1+te_2 & E_2^t=te_1 \\
\hline
{\bf C}(\alpha,\beta)  \to  {\bf A}_1(\alpha) &  
E_1^t=te_1+e_2 & E_2^t=t^2e_2\\
\hline
{\bf D}_1(\alpha,\beta)  \to  {\bf B}_2(\alpha) &  
E_1^t=te_1& E_2^t=e_2  \\
\hline
{\bf D}_1(\alpha,\beta)  \to  {\bf B}_2(1-\alpha+\beta) &  
E_1^t=te_2 & E_2^t=e_1-e_2  \\
\hline
{\bf D}_1(\alpha,\beta)  \to  {\bf D}_2(\beta, -\beta) &  
E_1^t=e_1  & E_2^t=te_2  \\
\hline
{\bf D}_2(\beta,\gamma)   \to  {\bf A}_3 &  
E_1^t=te_1+te_2 & E_2^t=t^2e_1+(\beta+\gamma)t^2e_2 \\
\hline
{\bf D}_3(\beta,\gamma)   \to  {\bf A}_2 &  
E_1^t=\frac{t}{1-\beta-\gamma}e_1- e_2 & E_2^t=te_2 \\
\hline
{\bf D}_3(\beta,\gamma)   \to  {\bf D}_2(\beta,\gamma) &  
E_1^t=e_1 & E_2^t=te_2 \\
\hline
{\bf D}_4(\alpha,\beta)   \to  {\bf A}_1(\beta) &  E_1^t=e_1+te_2 & E_2^t=-t^2e_2 \\
\hline
{\bf D}_4(\alpha,\beta)   \to  {\bf B}_2(\alpha) &  E_1^t=te_1 & E_2^t=e_2 \\
\hline
{\bf E}_2(\alpha,\beta,\gamma) \to  {\bf A}_1(\alpha) &  
E_1^t=te_1+e_2 & E_2^t=(1-\beta-\gamma)t^2e_1 \\
\hline
{\bf E}_3(\alpha,\delta,\epsilon) \to  {\bf B}_2\left(\frac{1-\delta-(1-\alpha)\epsilon}{1-\epsilon}\right) &
E_1^t=te_1 & E_2^t=\frac{\epsilon e_1-e_2}{\epsilon-1} \\
\hline
{\bf E}_4 \to  {\bf B}_3 &  
E_1^t=e_1-e_2 & E_2^t=te_2 \\
\hline
{\bf E}_4 \to  {\bf E}_5(\alpha) &  
E_1^t= \alpha e_1+(1-\alpha)e_2& E_2^t=(\alpha-t)e_1+(1-\alpha+t)e_2 \\
\hline
\end{longtable}

\begin{longtable}{|L{0.45\textwidth}|L{0.45\textwidth}|}
\hline
\multicolumn{2}{|l|}{\textbf{Table 3}}\\
\hline
\mbox{non-degenerations} & \mbox{separating sets}\\
\hline\hline
\endfirsthead
\hline
\multicolumn{2}{|l|}{\textbf{Table 3} (continued)}\\
\hline
\mbox{non-degenerations} & \mbox{separating sets}\\
\hline\hline
\endhead
\hline\endfoot
\endlastfoot
{\bf A}_1(\alpha)  \not\to  {\bf B}_3,\allowbreak  {\bf E}_5(\beta)\,\,(\beta\neq\alpha) &
\left\{\mu\left| \begin{array}{@{}L{0.34\textwidth}@{}}
c_{22}^1=c_{22}^2=c_{12}^1=c_{21}^1=0,\allowbreak \, c_{12}^2+c_{21}^2=c_{11}^1, \quad
(1-\alpha)c_{12}^2=\alpha c_{21}^2
\end{array} \right.\right\} \\
\hline
{\bf A}_2  \not\to  {\bf E}_5(\beta) &
\left\{\mu\left| \begin{array}{@{}L{0.34\textwidth}@{}}
c_{22}^1=c_{22}^2=c_{12}^1=c_{21}^1=c_{11}^1=0,\allowbreak \, c_{12}^2+c_{21}^2=0
\end{array} \right.\right\} \\
\hline
{\bf A}_4(\alpha)  \not\to   {\bf B}_2(\gamma),\allowbreak  {\bf D}_2(\beta,\gamma),\allowbreak  {\bf E}_5(\beta)  &
\left\{\mu\left| \begin{array}{@{}L{0.34\textwidth}@{}} 
c_{22}^1=c_{22}^2=c_{12}^1+c_{21}^1=0,\allowbreak  c_{12}^2 +c_{21}^2= c_{11}^1, \quad
\alpha^2 c_{12}^1 c_{11}^2=  (c_{11}^1)^2
\end{array} \right.\right\} \\
\hline
{\bf B}_1(\gamma)  \not\to  {\bf B}_2(\beta)\,\,(\beta\neq \gamma),\allowbreak  {\bf D}_2(\beta,\delta),\allowbreak   {\bf E}_5(\alpha)  &
\left\{\mu\left| \begin{array}{@{}L{0.40\textwidth}@{}} 
c_{22}^1=c_{22}^2=0,\allowbreak   
c_{12}^2+c_{21}^2=-c_{11}^1,\allowbreak   \\
c_{11}^2(c_{21}^1+c_{12}^1)=-(c_{11}^1)^2,\quad \gamma c_{12}^1 = (1-\gamma)c_{21}^1 \end{array}  						  \right.\right\} \\
\hline
{\bf B}_2(\gamma)  \not\to  {\bf B}_3,\allowbreak  {\bf E}_5(\beta) &
\left\{\mu\left| \begin{array}{@{}L{0.34\textwidth}@{}}
c_{11}^1=c_{12}^1=c_{21}^1=c_{22}^1=c_{22}^2=0,\\
(1-\gamma)c_{12}^2=\gamma c_{21}^2
\end{array} \right.\right\} \\
\hline
{\bf C}(\alpha,\beta)  \not\to 
{\bf B}_2(\gamma),\allowbreak {\bf B}_3,\allowbreak  {\bf D}_2(\gamma,\delta),\allowbreak  {\bf E}_5(\gamma)\,\,(\gamma \neq \alpha ) 
&
\left\{\mu\left| \begin{array}{@{}L{0.34\textwidth}@{}} 
c_{22}^1=0,\quad   c_{21}^2+c_{12}^2=c_{11}^1,\\
c_{21}^1=\alpha c_{22}^2,\quad c_{12}^1=(1-\alpha) c_{22}^2,\allowbreak  \\
(\alpha c_{21}^2-(1-\alpha) c_{12}^2)^2= \beta^2 c_{11}^2c_{22}^2  \end{array}\right.\right\}\\
\hline
{\bf D}_1(\alpha,\beta)  \not\to  
{\bf B}_2(\gamma)\,\,(\gamma \not\in \{ \alpha,1-\alpha+\beta\} ),\allowbreak {\bf B}_{3},\allowbreak  \\
\begin{array}{@{}L{0.43\textwidth}@{}} 
\qquad \qquad\qquad\quad   {\bf D}_{2}(\gamma,\delta)\,\,\big((\gamma,\delta) \neq (\beta,-\beta)\big),\allowbreak  {\bf E}_{5} (\gamma)
\end{array}&
\left\{\mu\left| \begin{array}{@{}L{0.34\textwidth}@{}}
                         c_{22}^1=c_{22}^2=0,\allowbreak   \alpha c_{12}^1=(1-\alpha) c_{21}^1,\allowbreak  \\ 
                         (\alpha-\beta)c_{12}^2-(1-\alpha+\beta) c_{21}^2=\beta c_{11}^1,\\
                         c_{11}^1(c_{12}^2+c_{21}^2)=c_{11}^2(c_{12}^1+c_{21}^1)
                         \end{array} 	  \right.\right\} \\
\hline
{\bf D}_2(\beta,\gamma)  \not\to  {\bf B}_3,\allowbreak  {\bf E}_5(\delta) &
\left\{\mu\left| \begin{array}{@{}L{0.34\textwidth}@{}}
c_{22}^1=c_{22}^2=c_{12}^1=c_{21}^1=0,\\
c_{12}^2=\beta c_{11}^1,\allowbreak \, c_{21}^2=\gamma c_{11}^1
\end{array} \right.\right\} \\
\hline
{\bf D}_3(\beta,\gamma)  \not\to
{\bf B}_2(\delta),\allowbreak  {\bf D}_2(\delta,\epsilon)\,\,\big( (\delta,\epsilon) \neq (\beta,\gamma) \big),\allowbreak  {\bf E}_{5}(\alpha)
&
\left\{\mu\left|  \begin{array}{@{}L{0.34\textwidth}@{}} 
c_{22}^1=c_{22}^2=c_{12}^1+c_{21}^1=0,\\
c_{12}^2+c_{21}^2=(\beta+\gamma)c_{11}^1,\\
(1-\beta-\gamma)(c_{12}^2-\beta c_{11}^1)c_{11}^1=c_{11}^2c_{12}^1
\end{array} \right.\right\} \\

\hline
{\bf D}_4(\alpha,\beta)  \not\to {\bf A}_1(\gamma)\,\,(\gamma\neq\beta),  {\bf A}_2,\, 
 {\bf B}_2(\gamma)\,\,(\gamma\neq\alpha),        &
\left\{\mu\left| \begin{array}{@{}L{0.4\textwidth}@{}}
c_{22}^1=c_{22}^2=0,\allowbreak \, \alpha c_{12}^1=(1-\alpha)c_{21}^1,\\
(1+\alpha-\beta)c_{12}^2+(1-\alpha-\beta)c_{11}^1=(1-\alpha+\beta)c_{21}^2,\\
(c_{11}^1-c_{12}^2-c_{21}^2)^2+4c_{11}^2(c_{12}^1+c_{21}^1)=0
\end{array} \right.\right\} \\
\qquad \qquad \qquad\qquad {\bf B}_3,\allowbreak  {\bf D}_2(\gamma,\delta),\allowbreak  {\bf E}_5(\gamma)\,\,(\gamma\neq\beta)  &\\
\hline
{\bf E}_1(\Gamma)  \not\to
  {\bf B}_2(\gamma),\allowbreak  {\bf B}_{3},\allowbreak  {\bf D}_{2}(\beta,\gamma)\,\,\big((\beta,\gamma) \not\in \mathcal{C}(\Gamma)\big),\allowbreak {\bf E}_5(\alpha)
 &
\mathcal{G}(\Gamma)\\

\hline
{\bf E}_2(\alpha,\beta,\gamma)  \not\to 
{\bf B}_2(\delta),\allowbreak  {\bf B}_{3},\allowbreak
\begin{array}{@{}L{0.43\textwidth}@{}}
   {\bf D}_2(\delta, \epsilon)\,\,\big((\delta, \epsilon)  \neq (\beta, \gamma)\big),  {\bf E}_5(\delta)\,\,(\delta \neq \alpha)
\end{array}
&
\mathcal{G}(1-\alpha,\beta,\alpha,\gamma) \\

\hline
{\bf E}_3(\alpha,\delta,\epsilon)  \not\to
{\bf B}_2(\gamma)\,\,\left(\gamma\not=\frac{1-\delta-(1-\alpha)\epsilon}{1-\epsilon}\right),\allowbreak   
\begin{array}{@{}L{0.43\textwidth}@{}}
\qquad\qquad {\bf B}_{3},  {\bf D}_2(\beta,\gamma)\,\,\big((\beta,\gamma)  \not\in \mathcal{C}(\alpha,\delta,\epsilon)\big),\allowbreak   {\bf E}_5(\gamma)
\end{array}
&
\mathcal{G}\left((1-\alpha)\epsilon,\frac{\delta}{\epsilon},\alpha\epsilon,\frac{1-\delta}{\epsilon}\right) \\
\hline
{\bf E}_4 \not\to
{\bf A}_3  &
\left\{\mu\left|  \begin{array}{@{}L{0.4\textwidth}@{}} 
c_{22}^1=c_{11}^2= 0,\allowbreak  c_{12}^1+c_{21}^1=c_{22}^2,\,  
c_{12}^2+c_{21}^2=c_{11}^1 
\end{array} \right.\right\} \\
\hline
\end{longtable}

\begin{longtable}{|L{0.15\textwidth}|L{0.77\textwidth}|}
\hline
\multicolumn{2}{|l|}{\textbf{Table 4}}\\
\hline
\endfirsthead
\hline
\multicolumn{2}{|l|}{\textbf{Table 4} (continued)}\\
\hline
\endhead
\hline\endfoot
\endlastfoot
{\bf A}_1(*) & {\bf A}_1(*),\allowbreak  {\bf A}_2,\allowbreak  {\bf A}_3,\allowbreak  {\bf B}_3,\allowbreak   {\bf E}_5(*),\allowbreak   {\bf k}^2 \\
\hline
{\bf A}_4(*)& {\bf A}_1(*),\allowbreak  {\bf A}_2,\allowbreak  {\bf A}_3,\allowbreak  {\bf A}_4(*),\allowbreak     {\bf B}_3,\allowbreak  {\bf E}_4,\allowbreak  {\bf E}_5(*),\allowbreak  {\bf k}^2 \\
\hline
{\bf B}_1(*)& {\bf A}_2,\allowbreak  {\bf A}_3,\allowbreak  {\bf A}_4',\allowbreak  {\bf B}_1(*),\allowbreak  {\bf B}_2(*),\allowbreak  {\bf B}_3,\allowbreak   {\bf k}^2 \\
\hline
{\bf B}_2(*) & {\bf A}_2,\allowbreak  {\bf A}_3,\allowbreak  {\bf B}_2(*),\allowbreak  {\bf B}_3,\allowbreak  {\bf k}^2  \\
\hline
{\bf C}(*) & {\bf A}_1(*),\allowbreak  {\bf A}_2,\allowbreak  {\bf A}_3,\allowbreak  {\bf A}_4(*),\allowbreak   {\bf B}_3,\allowbreak   {\bf C}(*),\allowbreak  {\bf E}_4,\allowbreak  {\bf E}_5(*),\allowbreak  {\bf k}^2  \\
\hline
{\bf D}_1(*) & {\bf A}_2,\allowbreak  {\bf A}_3,\allowbreak  {\bf A}_4',\allowbreak  {\bf B}_1(*),\allowbreak  {\bf B}_2(*),\allowbreak  {\bf B}_3,\allowbreak  {\bf D}_1(*),\allowbreak  {\bf D}_2'(*),\allowbreak  {\bf D}_3'(*),\allowbreak  {\bf k}^2  \\
\hline
  {\bf D}_2(*) & {\bf A}_1(*),\allowbreak  {\bf A}_2,\allowbreak  {\bf A}_3,\allowbreak  {\bf B}_2(*),\allowbreak  {\bf B}_3,\allowbreak  {\bf D}_2(*),\allowbreak  {\bf E}_5(*),\allowbreak  {\bf k}^2 \\
\hline
     {\bf D}_2'(*) & {\bf A}_2,\allowbreak  {\bf A}_3,\allowbreak  {\bf B}_3,\allowbreak  {\bf D}_2'(*),\allowbreak  {\bf k}^2 \\
\hline
  {\bf D}_3(*) & {\bf A}_1(*),\allowbreak  {\bf A}_2,\allowbreak  {\bf A}_3,\allowbreak  {\bf A}_4(*),\allowbreak  {\bf B}_1(*),\allowbreak   {\bf B}_2(*),\allowbreak {\bf B}_3,\allowbreak  {\bf D}_2(*),\allowbreak  {\bf D}_3(*),\allowbreak  {\bf E}_4,\allowbreak  {\bf E}_5(*),\allowbreak  {\bf k}^2 \\
\hline
     {\bf D}_3'(*) & {\bf A}_2,\allowbreak  {\bf A}_3,\allowbreak  {\bf A}_4',\allowbreak  {\bf B}_3,\allowbreak  {\bf D}_2'(*),\allowbreak  {\bf D}_3'(*),\allowbreak  {\bf k}^2 \\
\hline
{\bf D}_4(*) & {\bf A}_1(*),\allowbreak  {\bf A}_2,\allowbreak  {\bf A}_3,\allowbreak  {\bf A}_4(*),\allowbreak  {\bf B}_1(*),\allowbreak  {\bf B}_2(*),\allowbreak  {\bf B}_3,\allowbreak  {\bf D}_4(*),\allowbreak  {\bf E}_4,\allowbreak  {\bf E}_5(*),\allowbreak  {\bf k}^2 \\
\hline
  {\bf E}_1(*) & 
       \mathcal{A}_2  \\
\hline
{\bf E}_2(*) & {\bf A}_1(*),\allowbreak  {\bf A}_2,\allowbreak  {\bf A}_3,\allowbreak   {\bf A}_4(*),\allowbreak  {\bf B}_1(*),\allowbreak  {\bf B}_2(*),\allowbreak  {\bf B}_3,\allowbreak 
     {\bf C}(*),\allowbreak  {\bf D}_2(*),\allowbreak  {\bf D}_3(*),\allowbreak  {\bf D}_4(*),\allowbreak  {\bf E}_2(*),\allowbreak  {\bf E}_4,\allowbreak  {\bf E}_5(*),\allowbreak  {\bf k}^2  \\
\hline
{\bf E}_3(*) & {\bf A}_1(*),\allowbreak  {\bf A}_2,\allowbreak  {\bf A}_3,\allowbreak  {\bf A}_4(*),\allowbreak   {\bf B}_1(*),\allowbreak  {\bf B}_2(*),\allowbreak  {\bf B}_3,\allowbreak  {\bf D}_1(*),\allowbreak  {\bf D}_2(*),\allowbreak  {\bf D}_3(*),\allowbreak  {\bf D}_4(*),\allowbreak  {\bf E}_3(*),\allowbreak  {\bf E}_4,\allowbreak  {\bf E}_5(*),\allowbreak  {\bf k}^2 \\
\hline
{\bf E}_5(*) & {\bf B}_3,\allowbreak  {\bf E}_5(*),\allowbreak  {\bf k}^2  \\
\hline
\end{longtable}

\begin{longtable}{|L{0.27\textwidth}|L{0.36\textwidth}|L{0.27\textwidth}|}
\hline
\multicolumn{3}{|l|}{\textbf{Table 5}}\\
\hline
\mbox{degenerations} & \mbox{parametrized bases} & \mbox{parametrized indices}\\
\hline\hline
\endfirsthead
\hline
\multicolumn{3}{|l|}{\textbf{Table 5} (continued)}\\
\hline
\mbox{degenerations} & \mbox{parametrized bases} & \mbox{parametrized indices}\\
\hline\hline
\endhead
\hline\endfoot
\endlastfoot
{\bf A}_1(*) \to {\bf A}_2  & 
E_1^t=te_1,\allowbreak  E_2^t=t^2e_2 & \epsilon(t)= \frac{1}{t} \\
\hline
{\bf A}_4(*) \to {\bf A}_1(\alpha)  & 
E_1^t=te_1+(1-\alpha) e_2,\allowbreak  E_2^t=t^2e_2  & \epsilon(t)=\frac{1}{t}\\
\hline
{\bf A}_4(*) \to {\bf E}_4  &
E_1^t=te_1  ,\allowbreak  E_2^t=te_1+e_2  & \epsilon(t)=\frac{1}{t}  \\
\hline
{\bf B}_1(*) \to {\bf A}_4'  & 
E_1^t=-t^{-1}e_1+te_2  ,\allowbreak  E_2^t=-t^2e_2  & \epsilon(t)=\frac{1}{t^2}\\
\hline
{\bf B}_2(*) \to {\bf A}_2  & 
E_1^t=e_1+te_2,\allowbreak  E_2^t=-t^2e_2 & \epsilon(t)= \frac{1}{t} \\
\hline
{\bf C}(*) \to {\bf A}_4(\alpha) & 
E_1^t=te_1+\alpha e_2,\allowbreak  E_2^t=t^2e_2 & \epsilon(t)=\left( -\frac{1}{t^2}, \frac{\alpha(1+t^2)}{t^3} \right)\\
\hline
{\bf D}_1(*) \to {\bf B}_1(\alpha) & 
E_1^t=te_1,\allowbreak  E_2^t=e_2 & \epsilon(t)=\left( \alpha, \frac{1}{t} \right)\\
\hline
{\bf D}_1(*) \to {\bf D}_3'(\alpha) & 
E_1^t=e_1,\allowbreak  E_2^t=t e_2 & \epsilon(t)=\left( -\frac{1}{t}, \alpha \right)\\
\hline
{\bf D}_2(*) \to {\bf A}_1(\alpha)  & 
E_1^t=e_1+e_2,\allowbreak  E_2^t=te_2 & \epsilon(t)=\left( \alpha+t,  1- \alpha \right) \\
\hline
{\bf D}_2(*) \to {\bf B}_2(\alpha)  & 
E_1^t=e_2,\allowbreak  E_2^t=te_1 & \epsilon(t)=\left( \frac{\alpha}{t},  \frac{1- \alpha}{t} \right) \\
\hline
{\bf D}_2'(*) \to {\bf A}_2  & 
E_1^t=te_1-e_2,\allowbreak  E_2^t=te_2 & \epsilon(t)=\frac{1}{t}\\
\hline
{\bf D}_3(*) \to {\bf A}_4(\alpha)\,\,(\alpha\neq 0)  & 
E_1^t=\alpha e_1+\frac{1}{\alpha t}e_2  ,\allowbreak  E_2^t=e_2  & \epsilon(t)=\left(1+t+\frac{1}{\alpha^2t},  -\frac{1}{\alpha^2t}\right)\\
\hline
{\bf D}_3(*) \to {\bf B}_1(\alpha)  & 
E_1^t=-e_2  ,\allowbreak  E_2^t=te_1  & \epsilon(t)=\left(\frac{\alpha}{t},\frac{1-\alpha}{t}\right)\\
\hline
{\bf D}_3'(*) \to {\bf A}_4'  & 
E_1^t=te_1-\frac{1}{t}e_2 ,\allowbreak  E_2^t=e_2 & \epsilon(t)=-\frac{1}{t^2} \\
\hline
{\bf D}_4(*) \to {\bf A}_4(\alpha)\,\,(\alpha\neq0)  &
E_1^t=\alpha e_1+te_2,\allowbreak  E_2^t=-t^2e_2 & \epsilon(t)=\left(\frac{1}{t^2}, 1-\frac{1}{\alpha t}\right) \\
\hline
{\bf D}_4(*) \to {\bf A}_4'  &
E_1^t=-te_1+te_2,\allowbreak  E_2^t=-t^2e_2 & \epsilon(t)=\left(\frac{1}{t^2}, 1+\frac{1}{t^2}\right) \\
\hline
{\bf D}_4(*) \to {\bf B}_1(\alpha)  &
E_1^t=te_1-te_2,\allowbreak  E_2^t=e_2 & \epsilon(t)=\left(\alpha, \frac{1}{t}-1+\alpha\right) \\
\hline
{\bf E}_2(*) \to {\bf C}(\alpha,\beta)  & 
E_1^t=t^{-1}e_1-t^{-1}e_2 ,\allowbreak  E_2^t=e_2  & 
\epsilon(t)=\left(\alpha,\alpha+\beta t,1-\alpha-\beta t-t^2  \right)\\
\hline
{\bf E}_2(*) \to {\bf D}_3(\alpha,\beta)  & 
E_1^t=e_1,\allowbreak  E_2^t=te_2 & 
\epsilon(t)=\left(- \frac{1}{t},  \alpha, \beta \right)  \\
\hline
{\bf E}_2(*) \to {\bf D}_4(\alpha,\beta)  &
E_1^t=e_2,\allowbreak  E_2^t=te_1 & \epsilon(t)=\left(\beta, \frac{\alpha}{t}, \frac{1-\alpha}{t}\right) \\
\hline
{\bf E}_3(*) \to {\bf D}_1(\alpha,\beta)  & 
E_1^t=e_1,\allowbreak  E_2^t=te_2& 
\epsilon(t)=\left(\alpha,\frac{\beta}{t}  ,\frac{1}{t} \right)  \\
\hline
{\bf E}_3(*) \to {\bf D}_3(\alpha,\beta)\,\,(\alpha+\beta\neq 0)  & 
E_1^t=e_1 ,\allowbreak  E_2^t=-te_2  &
\epsilon(t)=\left(\frac{\alpha+\beta}{t},\frac{\alpha}{\alpha+\beta},\frac{1}{\alpha+\beta}\right)  \\
\hline
{\bf E}_3(*) \to {\bf D}_4(\alpha,\beta)  &
E_1^t=e_1,\allowbreak  E_2^t=\frac{(1+t)e_2-e_1}{t} & \epsilon(t)=\left(t(\alpha-1)+\beta, \frac{\beta}{1+t}, \frac{1}{1+t}\right) \\
\hline
{\bf E}_5(*) \to {\bf B}_3  & 
E_1^t=te_1,\allowbreak  E_2^t=e_2-e_1 & \epsilon(t)=\frac{1}{t} \\
\hline
\end{longtable}

\begin{longtable}{|L{0.46\textwidth}|L{0.45\textwidth}|}
\hline
\multicolumn{2}{|l|}{\textbf{Table 6}}\\
\hline
\mbox{non-degenerations} & \mbox{separating sets}\\
\hline\hline
\endfirsthead
\hline
\multicolumn{2}{|l|}{\textbf{Table 6} (continued)}\\
\hline
\mbox{non-degenerations} & \mbox{separating sets}\\
\hline\hline
\endhead
\hline\endfoot
\endlastfoot
{\bf A}_4(*) \not\to {\bf D}_4(\alpha,\beta) &
\left\{\mu\left| \begin{array}{@{}L{0.34\textwidth}@{}} c_{22}^1=c_{22}^2=c_{12}^1+c_{21}^1=0
\end{array} \right.\right\} \\
\hline
{\bf B}_1(*)  \not\to  {\bf D}_2'(\alpha)  &
\left\{\mu\left| \begin{array}{@{}L{0.34\textwidth}@{}}
c_{22}^1=c_{22}^2=c_{12}^1+c_{21}^1=c_{11}^1=0
\end{array}\right.\right\} \\
\hline
{\bf B}_1(*) \not\to {\bf D}_4(\alpha,\beta) &
\left\{\mu\left| \begin{array}{@{}L{0.34\textwidth}@{}} c_{22}^1=c_{22}^2=c_{12}^1+c_{21}^1=0
\end{array} \right.\right\} \\

\hline
{\bf C}(*) \not\to 
{\bf B}_2(\alpha),\allowbreak {\bf D}_2(\alpha,\beta)\,\,(\alpha+\beta\not=1),\allowbreak  {\bf D}_4(\alpha,\beta)
 & 
\left\{\mu\left| \begin{array}{@{}L{0.34\textwidth}@{}}
c_{22}^1=0,\allowbreak c_{12}^1+c_{21}^1=c_{22}^2,\allowbreak c_{12}^2+c_{21}^2=c_{11}^1
\end{array}\right.\right\}\\
\hline
{\bf D}_1(*) \not\to 
  {\bf A}_4(\alpha)\,\,(\alpha\not=0),\allowbreak  {\bf D}_2(\alpha,\beta)\,\,(\alpha+\beta\neq 0),\allowbreak 
 & 
\left\{\mu\left|  \begin{array}{@{}L{0.34\textwidth}@{}}
c_{22}^1=c_{12}^1+c_{21}^1=c_{11}^1=0
  \end{array}   \right.\right\} \\
\qquad\qquad\qquad\qquad{\bf D}_4(\alpha,\beta),\allowbreak {\bf E}_5(\alpha)&\\

\hline
{\bf D}_2(*) \not\to {\bf D}_4(\alpha,\beta) &
\left\{\mu\left| \begin{array}{@{}L{0.34\textwidth}@{}} c_{22}^1=c_{22}^2=c_{12}^1+c_{21}^1=0
\end{array} \right.\right\} \\
\hline
{\bf D}_2(*) \not\to {\bf A}_4(\alpha),\allowbreak  {\bf B}_1(\alpha),\allowbreak  {\bf D}_3(\alpha,\beta),\allowbreak  {\bf E}_4 &
\left\{\mu\left| \begin{array}{@{}L{0.34\textwidth}@{}} c_{22}^1=c_{22}^2=c_{12}^1=c_{21}^1=0
\end{array} \right.\right\} \\
\hline
{\bf D}_3(*) \not\to {\bf D}_4(\alpha,\beta) &
\left\{\mu\left| \begin{array}{@{}L{0.34\textwidth}@{}} c_{22}^1=c_{22}^2=c_{12}^1+c_{21}^1=0
\end{array} \right.\right\} \\
\hline
{\bf D}_3'(*) \not\to {\bf B}_2(\alpha)
& 
\left\{\mu\left|  \begin{array}{@{}L{0.34\textwidth}@{}} c_{22}^1=c_{22}^2=c_{12}^1+c_{21}^1=c_{12}^2+c_{21}^2=0
\end{array} \right.\right\} \\
\hline
{\bf D}_4(*) \not\to {\bf C}(\alpha,\beta),\allowbreak  {\bf D}_1(\alpha,\beta),\allowbreak  {\bf D}_2(\alpha,\beta),\allowbreak  {\bf D}_3(\alpha,\beta),  &
\left\{\mu\left| \begin{array}{@{}L{0.34\textwidth}@{}} c_{22}^1=c_{22}^2=0,\\ (c_{11}^1-c_{12}^2-c_{21}^2)^2+4c_{11}^2(c_{12}^1+c_{21}^1)=0
\end{array} \right.\right\} \\
\qquad\qquad \qquad \qquad  {\bf E}_1(\Gamma),\allowbreak  {\bf E}_2(\alpha,\beta,\gamma),\allowbreak  {\bf E}_3(\alpha,\beta,\gamma)& \\
\hline
{\bf E}_2(*) \not\to 
  {\bf D}_1(\alpha,\beta),\allowbreak  {\bf E}_1(\Gamma)\,\,(\Gamma\in\mathcal{V}),\allowbreak  {\bf E}_3(\alpha,\beta,\gamma)\,\,(\gamma\not=1)    & 
  \left\{\mu\left|  \begin{array}{@{}L{0.34\textwidth}@{}} c_{22}^1=0,\allowbreak  c_{12}^1+c_{21}^1=c_{22}^2 
\end{array} \right.\right\}
\\
\hline
{\bf E}_3(*) \not\to 
 {\bf C}(\alpha,\beta),\    {\bf E}_1(\Gamma)\,\,(\Gamma\in\mathcal{V}),&\\  
 \qquad\qquad\qquad{\bf E}_2(\alpha,\beta,\gamma)\,\,(\beta+\gamma\not=1)    & 
 \left\{\mu\left|  \begin{array}{@{}L{0.34\textwidth}@{}} c_{22}^1=c_{22}^2=0
\end{array} \right.\right\} 
\\
\hline
\end{longtable}

\endgroup
\section*{Appendix B. Figures.}
\addcontentsline{toc}{section}{Appendix B. Figures}

\begin{center}

{\bf Figure 1.}
\nopagebreak

\ 
\nopagebreak

\begin{adjustbox}{max width=\textwidth}\begin{tikzpicture}[->,>=stealth',shorten >=0.05cm,auto,node distance=0.94cm,
                    thick,main node/.style={rectangle,draw,fill=gray!10,rounded corners=1.5ex,font=\sffamily \scriptsize \bfseries },rigid node/.style={rectangle,draw,fill=black!20,rounded corners=1.5ex,font=\sffamily \scriptsize \bfseries },style={draw,font=\sffamily \scriptsize \bfseries }]
\node (0)   {$0$};
\node (0aa) [below  of=0]      {};
\node (0a1) [below  of=0aa]      {};
\node (0a) [below  of=0a1]      {};
\node (0a2) [below  of=0a]      {};
\node (0a3) [below  of=0a2]      {};
\node (0a4) [below  of=0a3]      {};
\node (0a5) [below  of=0a4]      {};
\node (0a6) [below  of=0a5]      {};
\node (0a7) [below  of=0a6]      {};
\node (1) [below  of=0a7]      {$1$};
\node (1aa) [below  of=1]      {};
\node (1a) [below  of=1aa]      {};
\node (2) [below  of=1a]      {$2$};
\node (2a) [below  of=2]      {};
\node (4) [below  of=2]      {$4$};

	\node            (Aa)  [right of =0]   {};                     

	\node[main node] (A4)  [right of =Aa]                       {${\bf A}_{4}(\alpha)$ };
	\node  (B4)  [right of =A4]   { };                     

	\node[main node] (A5)  [right of =B4]                       {${\bf B}_{1}(\gamma)$ };

    \node  (C4)  [right of =A5]   { };                     
	
	\node [main node]  (A55)  [right of =C4]                   {${\bf C}(\alpha,\beta)$ };                      
	
	\node  (B5)  [right of =A55]   { };                     
	
	\node[main node] (A8)  [right of =B5]                       {${\bf D}_{1}(\alpha,\beta)$ };

	\node  (B9)  [right of =A8]   { };                     
	
	\node[main node] (A10)  [right of =B9]                       {${\bf D}_{3}(\beta,\gamma )$ };
	\node (A19)  [right of =A10] {};
	\node[main node] (A18)  [right of =A19] {${\bf D}_4(\alpha,\beta)$};

	\node  (C12)  [right of =A18]   { };                     
	\node  (B12)  [right of =C12]   { };                     
	
	\node[main node] (A11)  [right of =C12]                       {${\bf E}_1(\Gamma)$};
	
	\node (A13)  [right of =A11]                       {};
	
	\node[main node] (A14)  [right of =A13]                       {${\bf E}_2(\alpha,\beta,\gamma)$};
	
	\node (A15)  [right of =A14]                       {};
	
	\node[main node] (A16)  [right of =A15]                       {${\bf E}_3(\alpha,\delta,\epsilon)$};
	
	\node (A17)  [right of =A16]                       {};
	
	\node[main node] (A12)  [right of =A17]                       {${\bf E}_4$};

	\node            (Aaa)  [right of =1]   {};                     
	\node  (Aaa21)  [right of =Aaa]   { };                     
	\node[main node] (A1)  [right of =Aaa21]                       {${\bf A}_{1}(\alpha)$ };
	\node  (B1)  [right of =A1]   { };                     

	\node[main node] (A2)  [right of =B1]                       {${\bf A}_{2}$ };
	\node  (B2)  [right of =A2]   { };                     

	\node[main node] (A6)  [right of =B2]                       {${\bf B}_{2}(\gamma)$ };
	\node  (B6)  [right of =A6]   { };                     
	\node  (B8)  [right of =B6]   { };                     
	\node  (B10)  [right of =B8]   { };
	\node  (B12)  [right of =B10]   { };
	\node  (B1222)  [right of =B12]   { };

	\node[main node] (A9)  [right of =B1222]                       {${\bf D}_{2}(\beta,\gamma)$ };
	\node  (B2)  [right of =A2]   { };

	\node  (21)  [right of =2]   { };                     
	\node  (22)  [right of =21]    { };                     
	\node  (221)  [right of =22]   { };                     
	\node  (23)  [right of =221]   { };                     
	\node  (24)  [right of =23]   { };                     
	
	\node[main node] (A3)  [right of =24]                       {${\bf A}_{3}$ };
	
	\node  (25)  [right of =A3]   { };                     
	\node  (251)  [right of =25]   { }; 
	\node  (2511)  [right of =251]   { };
	\node  (25111)  [right of =2511]   { };
	\node  (251111)  [right of =25111]   { };
	
	\node[main node] (A7)  [right of =251111]                       {${\bf B}_{3}$ };

	\node  (26)  [right of =A7]   { };                     
	\node  (261)  [right of =26]   { };                     

	\node[main node] (A11a1a)  [right of =261]                       {${\bf E}_5(\alpha)$ };

\node (1160)[right of=4]{};
\node (1161)[right of=1160]{};
\node (1162)[right of=1161]{};
\node (1163)[right of=1162]{};
\node (1164)[right of=1163]{};
\node (1165)[right of=1164]{};
\node (1166)[right of=1165]{};
\node (1167)[right of=1166]{};

	\node[main node] (C) [ right  of=1167]       {${\bf k}^2$};
 
\path[every node/.style={font=\sffamily\small}]


  (A1)   edge  (A3) 
  (A2)   edge  (A3) 
  (A2)   edge  (A7) 
  (A4)   edge  (A2) 
  (A55)  edge  node {} (A1) 
  (A5)   edge  (A2) 
  (A6)   edge  (A3) 
  (A8)   edge [bend right=0] node[sloped,pos=0.29,fill=white,inner sep=1pt]{\scriptsize $\gamma\in\{\alpha,1-\alpha+\beta\}$ }    node {} (A6) 
  (A8)   edge  [bend right=0]  node[sloped,pos=0.68,fill=white,inner sep=1pt]{\scriptsize $\beta+\gamma=0$ }   node {} (A9) 
  (A5)   edge  (A6) 
  (A10)   edge [bend right=-8]     node {} (A2) 

  (A12)   edge [bend right=0]     node {} (A7)
  
  (A12)   edge [bend right=0]     node {} (A11a1a) 


   (A14)  edge [bend right=0] node {} (A9) 

   (A14) edge [bend right=0]
   node {} (A1) 

  (A9)   edge  (A3) 
  (A10)  edge  (A9) 
  (A11)  edge [bend right=0] node[sloped,pos=0.56,fill=white,inner sep=1pt]{\scriptsize  $(\beta,\gamma)\in\mathcal{C}(\Gamma)$  }   node {} (A9)
  
  (A16)  edge [bend right=0] node[sloped,pos=0.50,fill=white,inner sep=1pt]{\scriptsize  $(\beta,\gamma)\in\mathcal{C}(\alpha,\delta,\epsilon)$  }   node {} (A9)
  
  (A16)  edge [bend right=7] node[sloped,pos=0.73,fill=white,inner sep=1pt]{\scriptsize  $\gamma=\frac{1-\delta-(1-\alpha)\epsilon}{1-\epsilon}$   }   node {} (A6)

  (A1)   edge  (A11a1a) 
  (A18)  edge[newedge] node[sloped,pos=0.20,fill=white,inner sep=1pt]{\scriptsize ${\bf A}_1(\beta)$} (A1)
  (A18)  edge[newedge] node[sloped,pos=0.63,fill=white,inner sep=1pt]{\scriptsize ${\bf B}_2(\alpha)$} (A6)
  
  (A3)   edge  (C) 
  (A7)   edge  (C) 
(A11a1a) edge  (C)  
 
 ;

\end{tikzpicture}\end{adjustbox}


\end{center}

\

\begin{center}

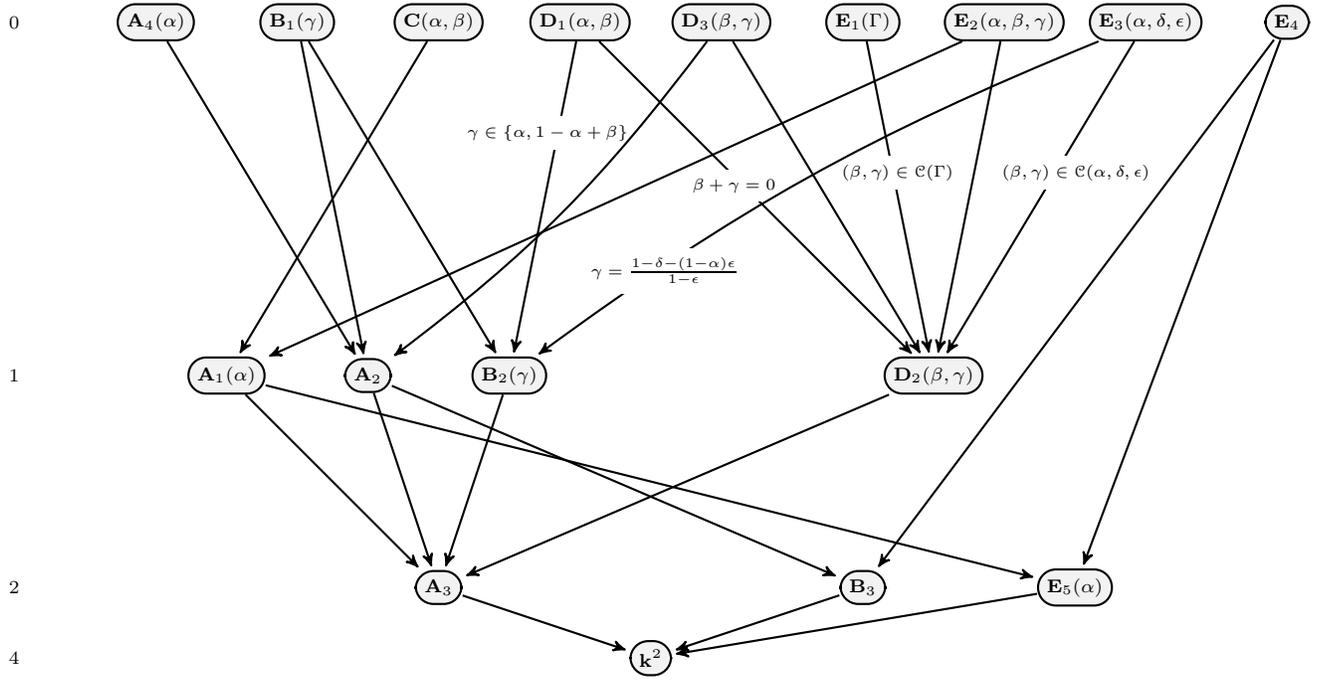
{\bf Figure 2.}
\nopagebreak

\ 
\nopagebreak

\begin{adjustbox}{max width=\textwidth}\begin{tikzpicture}[-,draw=gray!50,node distance=0.93cm,
                   ultra thick,main node/.style={rectangle, fill=gray!50,font=\sffamily \scriptsize \bfseries },style={draw,font=\sffamily \scriptsize \bfseries }]

\node (8)   {$8$};
\node (8r) [right  of=8] {};
\node (7) [right  of=8]      {$7$};
\node (7r) [right  of=7] {};
\node (6) [right  of=7r]      {$6$};
\node (6r) [right  of=6] {};
\node (6rr) [right  of=6r] {};
\node (6rrr) [right  of=6rr] {};
\node (5) [right  of=6r]      {$5$};
\node (5r) [right  of=5] {};
\node (5rr) [right  of=5r] {};
\node (5rrr) [right  of=5rr] {};
\node (4) [right  of=5r]      {$4$};
\node (4r) [right  of=4] {};
\node (3) [right  of=4r]      {$3$};
\node (3r) [right  of=3] {};
\node (2) [right  of=3]      {$2$};
\node (2r) [right  of=2] {};
\node (0) [right  of=2]      {$0$};

\node (8b) [below of =8] {};
\node (8bb) [below of =8b] {};

	\node[main node] (A11)  [below of =8bb]                       {$\overline{O\big({\bf E}_1(*)\big)}$};
	
\node (7b) [below of =7] {};
	
	\node[main node] (A14)  [below of =7b]                       {$\overline{O\big({\bf E}_3(*)\big)}$};
	
\node (A14b) [below of =A14] {};
	
	\node[main node] (A16)  [below of =A14b]                       {$\overline{O\big({\bf E}_2(*)\big)}$};

\node (6b) [below of =6] {};

	\node[main node] (A8)  [below of =6b]                       {$\overline{O\big({\bf D}_{1}(*)\big)}$ };

	\node[main node] (A10)  [below of =A8]                       {$\overline{O\big({\bf D}_{3}(*)\big)}$ };
	\node[main node] (D4n)  [below of =A10] {$\overline{O\big({\bf D}_{4}(*)\big)}$};
	
	\node [main node]  (A55)  [below of =D4n]                   {$\overline{O\big({\bf C}(*)\big)}$ };

	\node[main node] (A5)  [below of =5]                       {$\overline{O\big({\bf B}_{1}(*)\big)}$ };
	
	\node[main node] (A9)  [below of =A5]                       {$\overline{O\big({\bf D}_{2}(*)\big)}$ };
	
\node (A5b) [below of =A9] {};

	\node[main node] (A10s)  [below of =A5b]                       {$\overline{O\big({\bf D}'_{3}(*)\big)}$ };
	
	\node[main node] (A4)  [below of =A10s]                       {$\overline{O\big({\bf A}_{4}(*)\big)}$ };

	\node[main node] (A6)  [below of =4]                       {$\overline{O\big({\bf B}_{2}(*)\big)}$ };
	
	\node[main node] (A9s)  [below of =A6]                       {$\overline{O\big({\bf D}'_{2}(*)\big)}$ };
	
	\node[main node] (A4s)  [below of =A9s]                       {$\overline{O(A_4')}$ };
	
	\node[main node] (A1)  [below of =A4s]                       {$\overline{O\big({\bf A}_{1}(*)\big)}$ };

	\node[main node] (A12)  [below of =A1]                       {$\overline{O({\bf E}_{4})}$ };
	
\node (3b) [below of =3] {};
	
	\node[main node] (A2)  [below of =3b]                       {$\overline{O({\bf A}_{2})}$ };
	
\node (A2b) [below of =A2] {};
	
	\node[main node] (E5)  [below of =A2b]                       {$\overline{O\big({\bf E}_5(*)\big)}$ };
	
\node (2b) [below of =2] {};
	
	\node[main node] (A3)  [below of =2]                       {$\overline{O({\bf A}_{3})}$ };
	
\node (A3b) [below of =A3] {};
	
	\node[main node] (A7)  [below of =A3b]                       {$\overline{O({\bf B}_{3})}$ };
	
\node (0b) [below of =0] {};
\node (0bb) [below of =0b] {};
	
	\node[main node] (C) [ below  of =0b]       {${\bf k}^2$};
 


\path (C) edge   (A3)


 (C) edge   (A7) 

 (A3) edge   (A2) 

 (A7) edge   (A2) 

 (A7) edge   (E5) 

(A2) edge   (A1)

(A2) edge   (A4s)

(A2) edge   (A6)

(A2) edge   (A9s)





 (E5) edge   (A12) 

 (E5) edge   (A1) 

 (A1) edge   (A4) 

 (A1) edge   (A9) 

 (A4s) edge   (A4) 

 (A4s) edge   (A5) 

 (A4s) edge   (A10s) 

 (A6) edge   (A5) 

 (A6) edge   (A9) 

 (A9s) edge   (A9) 

 (A9s) edge   (A10s) 

 (A12) edge   (A4) 

 (A4) edge   (A10) 

 (A4) edge   (A55) 

 (A5) edge   (A10) 

 (A9) edge   (A10) 

 (A10s) edge   (A10) 

 (A5) edge   (A8) 

 (A10s) edge   (A8) 

 (A8) edge   (A14) 

 (A10) edge   (A14) 

 (A10) edge   (A16) 

 (A55) edge   (A16) 

 (A4) edge[newedge]   (D4n) 

 (A5) edge[newedge]   (D4n) 

 (D4n) edge[newedge]   (A14) 

 (D4n) edge[newedge]   (A16) 

 (A14) edge   (A11) 

 (A16) edge   (A11) ;

\end{tikzpicture}\end{adjustbox}

\end{center}

\begin{center}
{
{\bf Figure 3.}
\nopagebreak

\ 
\nopagebreak

\begin{adjustbox}{max width=\textwidth}\begin{tikzpicture}[->,>=stealth',shorten >=0.05cm,auto,node distance=1cm,
                    thick,main node/.style={rectangle,draw,fill=gray!10,rounded corners=1.5ex,font=\sffamily \scriptsize \bfseries },rigid node/.style={rectangle,draw,fill=black!20,rounded corners=1.5ex,font=\sffamily \scriptsize \bfseries },style={draw,font=\sffamily \scriptsize \bfseries }]

\node (0)   {$0$};
\node (0aa) [below  of=0]      {};

\node (0b1) [below  of=0aa]      {};
\node (0b2) [below  of=0b1]      {};
\node (0b3) [below  of=0b2]      {};
\node (0b4) [below  of=0aa]      {};
\node (1) [below  of=0b4]      {$1$};
\node (1aa) [below  of=1]      {};

\node (2) [below  of=1aa]      {$2$};

\node (4) [below  of=2]      {$4$};

\node            (00)  [right of =0]   {};                     
\node            (01)  [right of =00]   {};                     
\node            (02)  [right of =01]   {};                     

\node[main node] (4c)  [right of =00]                       {${\bf C}(\frac{1}{2},0)$ };

\node            (021)  [right of =4c]   {};

\node[main node] (C1)  [right of =021]                       {${\bf D}_1(\frac{1}{2},0)$ };

\node            (10)  [right of =C1]   {};                     

\node[main node] (D4c) [right of =10] {${\bf D}_4(\frac{1}{2},\frac{1}{2})$};
\node (10d) [right of =D4c] {};
\node[main node] (C8)  [right of =10d]                       {${\bf E}_2^c(\alpha)$ };

\node            (11)  [right of =C8]   {};

\node[main node] (C9)  [right of =11]                       {${\bf E}_3^c(\beta)$ };

\node            (03)  [right of =C9]   {};
\node            (051)  [right of =03]   {};
\node            (052)  [right of =051]   {};                     

\node[main node] (C7)  [right of =051]                       {${\bf E}_1^c(\beta,\gamma)$ };

\node            (12)  [right of =1]   {};                     
\node            (12q)  [right of =12]   {};                     
\node            (12qq)  [right of =12q]   {};                     

	\node[main node] (A1)  [right of =12qq]                       {${\bf A}_1(\frac{1}{2})$ };
\node            (16)  [right of =A1]   {};                     

	\node[main node] (B2)  [right of =16]                       {${\bf B}_2(\frac{1}{2})$ };
\node            (17)  [right of =B2]   {};                     

	\node[main node] (C2)  [right of =17]                       {${\bf D}_2^c(\alpha)$ };

\node            (20)  [right of =2]   {};                     
\node            (21)  [right of =20]   {};                     
\node            (22)  [right of =21]   {};                     

\node[main node] (C212)  [right of =22]                       {${\bf E}_5(\alpha)$ };

\node            (25)  [right of =C212]   {};                     

	\node[main node] (A3)  [right of =25]                       {${\bf A}_{3}$ };

\node            (26)  [right of =A3]   {};                     

	\node[main node] (F)  [right of =26]                       {${\bf B}_3$ };

\node (q1)[right of=C212]{};
\node (q2)[right of=q1]{};
\node (q21)[right of=q2]{};
\node (q31)[right of=q21]{};

\node (q3)[right of=q31]   {};

\node (q4)[below  of=q3]   {};

\node (1160)[right of=4]{};
\node (1161)[right of=1160]{};
\node (1162)[right of=1161]{};
\node (1163)[right of=1162]{};
\node (1164)[right of=1163]{};
\node (1165)[right of=1164]{};

	\node[main node] (C) [ right  of=1164]       {${\bf k}^2$};
 
\path[every node/.style={font=\sffamily\small}]


  (C1)   edge  (B2) 
  (C1)   edge node[sloped,pos=0.62,fill=white,inner sep=1pt]{\scriptsize $\alpha=0$ }   node {} (C2) 
  


  (4c)   edge  (A1)

  (A1)   edge  (A3) 
  (A1)   edge node[left,pos=0.5,inner sep=2pt]{\scriptsize $\alpha=1/2$ }    (C212) 

  (B2)   edge  (A3) 

  (C2)   edge  (A3) 

  (C8)   edge  (A1) 
  (C9)   edge  (B2) 
  (D4c)  edge[newedge]  (A1) 
  (D4c)  edge[newedge]  (B2)

  (C7)   edge node[sloped,pos=0.5,fill=white,inner sep=1pt]{\scriptsize $\alpha\in\{\beta,\gamma,\frac{\beta+\gamma-1}{4\beta\gamma-1}\}$ }    (C2) 
  (C8)   edge  (C2) 
  (C9)   edge node[sloped,pos=0.65,fill=white,inner sep=1pt]{\scriptsize $\alpha\in\{\frac{\beta}{2},\frac{1}{2\beta}\}$ }     (C2)

  (A3)   edge  (C) 
  (F)   edge  (C) 

  (C212)   edge  (C) ;
\end{tikzpicture}\end{adjustbox}
}
\end{center}

\

\begin{center}

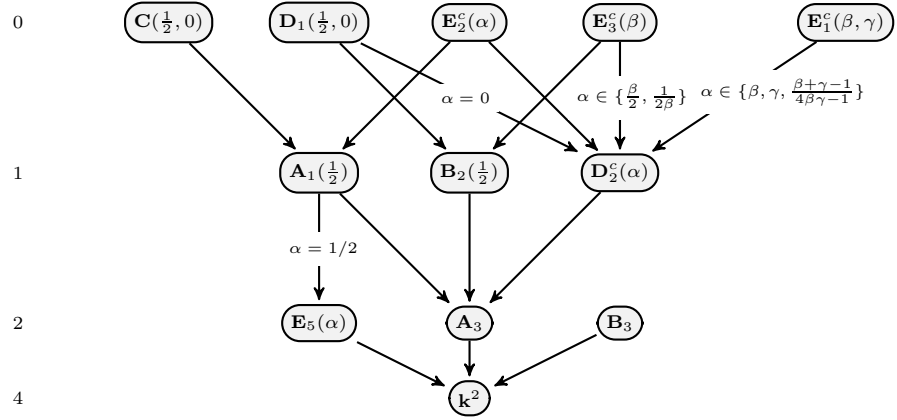
{\bf Figure 4.}
\nopagebreak

\ 
\nopagebreak

\begin{adjustbox}{max width=\textwidth}\begin{tikzpicture}[-,draw=gray!50,node distance=0.93cm,
                   ultra thick,main node/.style={rectangle, fill=gray!50,font=\sffamily \scriptsize \bfseries },style={draw,font=\sffamily \scriptsize \bfseries }]

\node (6) {$6$};
\node (6r) [right  of=6] {};
\node (6rr) [right  of=6r] {};
\node (5) [right  of=6r]      {$5$};
\node (5r) [right  of=5] {};
\node (5rr) [right  of=5r] {};
\node (4) [right  of=5r]      {$4$};
\node (4r) [right  of=4] {};
\node (4rr) [right  of=4r] {};

\node (3) [right  of=4r]      {$3$};
\node (3r) [right  of=3] {};
\node (3rr) [right  of=3r] {};
\node (2) [right  of=3r]      {$2$};
\node (2r) [right  of=2] {};
\node (2rr) [right  of=2r] {};

\node (0) [right  of=2r]      {$0$};

\node (6b) [below of =6] {};
\node (6bb) [below of =6b] {};

	\node[main node] (C7)  [below of =6b]                       {$\overline{O\big({\bf E}_1^c(*)\big)}$};

	\node[main node] (C9)  [below of =5]                       {$\overline{O\big({\bf E}_3^c(*)\big)}$};

\node (5b) [below of =C9] {};

	\node[main node] (C8)  [below of =5b]                       {$\overline{O\big({\bf E}_2^c(*)\big)}$};

	\node[main node] (C5)  [below of =4]                       {$\overline{O\big({\bf D}_1(\frac{1}{2},0)\big)}$};
	\node[main node] (C6)  [below of =C5]                       {$\overline{O\big({\bf D}_2^c(*)\big)}$};
	\node[main node] (D4c)  [below of =C6] {$\overline{O\big({\bf D}_4(\frac{1}{2},\frac{1}{2})\big)}$};
	\node[main node] (C4)  [below of =D4c]                       {$\overline{O\big({\bf C}(\frac{1}{2},0)\big)}$};

	\node[main node] (C3)  [below of =3]                       {$\overline{O\big({\bf B}_2(\frac{1}{2})\big)}$};
	\node[main node] (C60)  [below of =C3]                       {$\overline{O\big({\bf D}_2^c(0)\big)}$};
	\node[main node] (C1)  [below of =C60]                       {$\overline{O\big({\bf A}_1(\frac{1}{2})\big)}$};

\node (2b) [below of =2] {};

	\node[main node] (C2)  [below of =2b]                       {$\overline{O\big({\bf A}_3\big)}$};

    \node[main node] (C10)  [below of =C2]                       {$\overline{O\big({\bf E}_5(\frac{1}{2})\big)}$};
    \node[main node] (AA)  [below of =C10]                       {$\overline{O\big({\bf B}_3\big)}$};
    \node[main node] (R1)  [below of =C1]                       {$\overline{O\big({\bf E}_5(*)\big)}$};

\node (0b) [below of =0] {};
\node (0bb) [below of =0b] {};
	
	\node[main node] (C) [ below  of =0bb]       {${\bf k}^2$};
 


\path 


 (C7) edge   (C8) 
 (C7) edge   (C9)
 
 (C8) edge   (C4) 
 (C8) edge   (C6)
 
 (C9) edge   (C6)
 (C9) edge   (C5)
 
 (C5) edge   (C3)
 (C5) edge   (C60)
 
 (C6) edge   (C60)
 (C6) edge   (C3)
 (C6) edge   (C1)

 (C4) edge   (C1) 
 (C8) edge[newedge]   (D4c) 
 (C9) edge[newedge]   (D4c) 
 (D4c) edge[newedge]   (C3) 
 (D4c) edge[newedge]   (C1) 
 
 (C3) edge   (C2)
 (C60) edge   (C2)
 (C1) edge   (C10)
 (C1) edge   (C2) 

 (C2) edge   (C) 
 (C10) edge   (C) 
 (R1) edge   (C10) 
 (R1) edge   (AA) 
 (AA) edge   (C)

 ;
\end{tikzpicture}\end{adjustbox}

\end{center}

\begin{center}

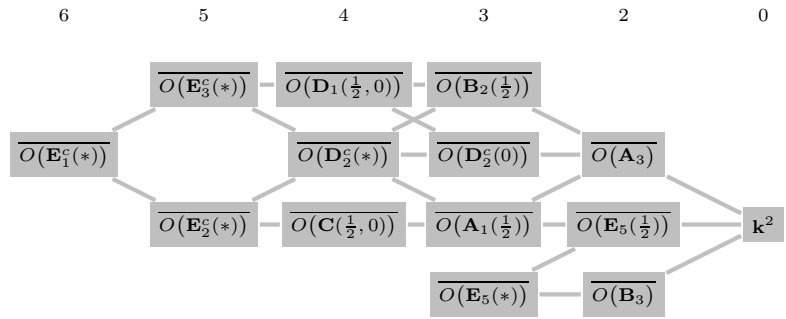
{\bf Figure 5.}
\nopagebreak

\ 
\nopagebreak

\begin{adjustbox}{max width=\textwidth}\begin{tikzpicture}[->,>=stealth',shorten >=0.05cm,auto,node distance=1cm,
                    thick,main node/.style={rectangle,draw,fill=gray!10,rounded corners=1.5ex,font=\sffamily \scriptsize \bfseries },rigid node/.style={rectangle,draw,fill=black!20,rounded corners=1.5ex,font=\sffamily \scriptsize \bfseries },style={draw,font=\sffamily \scriptsize \bfseries }]

\node (0)   {$0$};
\node (0a) [below  of=0]      {};
\node (0a1) [below  of=0aa]      {};
\node (0a2) [below  of=0a1]      {};
\node (1) [below  of=0]      {$1$};

\node (1a) [below  of=1]      {};
\node (1aa) [below  of=1a]      {};
\node (2) [below  of=1]      {$2$};
\node (2a) [below  of=2]      {};
\node (4) [below  of=2]      {$4$};

\node            (00)  [right of =0]   {};                     
\node            (01)  [right of =00]   {};                     
\node            (02)  [right of =01]   {};                     
\node            (03)  [right of =02]   {};                     

\node[main node] (L7)  [right of =03]                       {${\bf D}_1(0,0)$};
\node            (x0)  [right of =L7]   {};                     

\node            (x00)  [right of =x0]   {};                     

\node  (L9)  [right of =x00]                       { };
\node            (x1)  [right of =L9]   {};                     

\node  (L6)  [right of =x1]                       { };
\node            (x11)  [right of =L6]   {};

\node[main node] (L12)  [right of =x00]                     {${\bf E}_1(0,0,0,0)$};
\node            (x2)  [right of =L12]   {};                     
\node            (x3)  [right of =x2]   {};                     

\node (L13)  [right of =x2]                     { };

\node            (10)  [right of =1]   {};                     
\node            (11)  [right of =10]   {};                     
\node            (12)  [right of =11]   {};

\node[main node] (L3)  [right of =11]                       {${\bf B}_2(0)$};
\node            (y1)  [right of =L3]   {};                     

\node[main node] (L4)  [right of =y1]                       {${\bf B}_2(1)$};
\node            (y2)  [right of =L4]   {};                     

\node  (L1)  [right of =y2]                       { };
\node            (y3)  [right of =L1]   {};                     

\node[main node] (L10)  [right of =y3]                       {${\bf D}_2(1,1)$};
\node            (y4)  [right of =L10]   {};                     

\node[main node] (L11)  [right of =y2]                       {${\bf D}_2(0,0)$};

\node            (20)  [right of =2]   {};                     
\node            (21)  [right of =20]   {};                     
\node            (22)  [right of =21]   {};                     
\node            (z1)  [right of =22]   {};

\node  (L5)  [right of =z1]                       { };
\node            (z2)  [right of =L5]   {};                     

\node[main node] (L2)  [right of =L5]                       {${\bf A}_3$};
\node            (z3)  [right of =L2]   {};

\node (1160)[right of=4]{};
\node (1161)[right of=1160]{};
\node (1162)[right of=1161]{};
\node (1163)[right of=1162]{};
\node (1164)[right of=1163]{};
\node (1165)[right of=1164]{};

	\node[main node] (C) [ right  of=1164]       {${\bf k}^2$};
 
\path[every node/.style={font=\sffamily\small}]

(L3)  edge  (L2)

(L4)  edge  (L2)

(L7)  edge  (L3)
(L7)  edge  (L4)
(L7)    edge    (L11)


(L10)  edge  (L2)
(L11)  edge  (L2)

(L12)  edge    (L10)
(L12)  edge    (L11)

(L2)  edge  (C)
 ;
         
\end{tikzpicture}\end{adjustbox}
\end{center}

\section*{Appendix C. Differences from the published version}
\addcontentsline{toc}{section}{Appendix C. Differences from the published version}

\emph{Theorem~\ref{alg}: statement changed.} The family of algebras ${\bf D}_4(\alpha,\beta)$, missing in the
published classification, has been added. Accordingly, Theorems~\ref{graf2} and~\ref{inf2}, the tables and the
figures have been completed: the missing degenerations have been constructed and the missing non-degenerations have
been proved.
 
\noindent 
\emph{Theorem~\ref{inf2}: proof completed.} The non-degenerations ${\bf D}_2(*)\not\to{\bf A}_4(\alpha),{\bf B}_1(\alpha),{\bf D}_3(\alpha,\beta),{\bf E}_4$ were missing; the new row of Table~6 and the corresponding paragraph of the proof have been added.

\end{document}